\documentclass[11pt]{amsart}
\usepackage{amsmath, amssymb, latexsym}
\usepackage{mathrsfs}
\usepackage[all, knot]{xy}
\usepackage{bbm}
\usepackage{color, url}
\xyoption{arc}

\numberwithin{equation}{section}

\theoremstyle{plain}
\newtheorem{theorem}{Theorem}[section] 
\newtheorem{lemma}[theorem]{Lemma}
\newtheorem{proposition}[theorem]{Proposition}
\newtheorem{corollary}[theorem]{Corollary}

\theoremstyle{definition}
\newtheorem{definition}[theorem]{Definition}
\newtheorem{example}[theorem]{Example}

\theoremstyle{remark}
\newtheorem{remark}[theorem]{Remark}

\newcommand{\vs}{\vspace}
\newcommand{\hs}{\hspace}
\newcommand{\C}{\mathbb{C}}
\newcommand{\B}{\mathbb{B}}
\newcommand{\Q}{\mathbb{Q}}
\newcommand{\Z}{\mathbb{Z}}
\newcommand{\N}{\mathbb{N}}
\newcommand{\R}{\mathcal{R}}
\newcommand{\h}{\mathfrak{h}}

\newcommand{\PP}{\mathbf{P}}
\newcommand{\QQ}{\mathbf{Q}}
\newcommand{\KK}{\mathbf{K}}
\newcommand{\LL}{\mathcal{L}}
\newcommand{\Pp}{\mathcal{P}}
\newcommand{\A}{\mathcal A}
\newcommand{\EE}{\mathbf E}
\newcommand{\F}{\mathcal F}
\newcommand{\OO}{\mathcal O}
\newcommand{\FF}{\mathscr F}
\newcommand{\NN}{\mathscr N}
\newcommand{\E}{\mathcal E}
\newcommand{\PPP}{\mathscr{P}}
\newcommand{\ra}{\rightarrow}
\newcommand{\ma}{\mapsto}
\newcommand{\co}{\cong}
\newcommand{\ti}{\times}
\newcommand{\ot}{\otimes}
\newcommand{\pl}{\partial}
\newcommand{\op}{\bigoplus}
\newcommand{\Ga}{\Gamma}
\newcommand{\ga}{\gamma}
\newcommand{\La}{\Lambda}
\newcommand{\la}{\lambda}
\newcommand{\De}{\Delta}
\newcommand{\de}{\delta}
\newcommand{\om}{\omega}
\newcommand{\al}{\alpha}
\newcommand{\si}{\sigma}
\newcommand{\bt}{\beta}
\newcommand{\vp}{\varphi}
\newcommand{\ep}{\epsilon}
\newcommand{\vep}{\varepsilon}
\newcommand{\wi}{\widetilde}
\newcommand{\ov}{\overset}
\newcommand{\ovv}{\overline}
\newcommand{\se}{\subset}
\newcommand{\cs}{\cdots}
\newcommand{\ds}{\dots}
\newcommand{\ac}{\cdot}
\newcommand{\x}{\hs{0.5mm}\cdot\hs{0.5mm}}

\newcommand{\bu}{\bullet}
\newcommand{\tx}{\text}
\newcommand{\xs}{\xrightarrow{\sim}}
\newcommand{\xr}{\xrightarrow}
\newcommand{\s}{\displaystyle\sum}
\newcommand{\p}{\displaystyle\prod}
\newcommand{\f}{\displaystyle\frac}

\newcommand{\Hom}{\tx{Hom}}
\newcommand{\HOM}{\tx{HOM}}
\newcommand{\End}{\tx{End}}
\newcommand{\zz}{{-1}}
\newcommand{\ad}{\text{ad}}
\newcommand{\wt}{\text{wt}}
\newcommand{\Ind}{\text{Ind}}
\newcommand{\Res}{\text{Res}}

\newcommand{\Dim}{\mathbf{Dim}\ }
\newcommand{\Mod}{\text{-Mod}}
\newcommand{\pMod}{\text{-pMod}}
\newcommand{\fMod}{\text{-fMod}}
\newcommand{\Ext}{\text{Ext}}
\newcommand{\Spr}{\text{Spr}}
\newcommand{\Ch}{\text{Ch}}
\newcommand{\hd}{\text{hd}}

\newcommand{\Seq}{\text{Seq}}
\newcommand{\ii}{\textit{\textbf{i}}}
\newcommand{\jj}{\textit{\textbf{j}}}
\newcommand{\kk}{\textit{\textbf{k}}}
\newcommand{\II}{\mathbbm{I}}
\newcommand{\ZZ}{\mathscr{Z}}

\title[Quantum groups of Borcherds-Cartan and KLR algebras]
{Quantum groups of Borcherds-Cartan type and Khovanov-Lauda-Rouquier algebras}

\author[Seok-Jin Kang]{Seok-Jin Kang}
\address{Korea Research Institute of Arts and Mathematics,
Asan-si, Chungcheongnam-do, 31551, Korea}
\email{stefano.kang47@gmail.com}

\author[Young Rock Kim]{Young Rock Kim}
\address{Graduate School of Education, Hankuk University of Foreign Studies, Seoul, 02450,  Korea}
\email{rocky777@hufs.ac.kr}

\author[Bolun Tong]{Bolun Tong}
\address{Graduate School of Education, Hankuk University of Foreign Studies, Seoul, 02450, Korea}
\email{tbl\_2018@163.com}

\keywords{Categorification, KLR algebra, quiver with loops}

\subjclass[2010] {17B37, 17B67, 16G20}

\begin{document}

\maketitle

\begin{abstract}
We categorify a class of quantum groups associated with quivers, possibly with loops, by constructing the corresponding Khovanov-Lauda-Rouquier algebras (KLR) algebras $R$. We prove that the indecomposable projective $R$-modules realize the canonical basis of the negative part $U^-$ of the quantum group. Moreover, for $\La \in P^+$, the cyclotomic KLR algebra $R^\La$ provide a categorification of the irreducible highest weight $U$-module $V(\La)$.
\end{abstract}

\section*{Introduction}

{\it Khovanov-Lauda-Rouquier algebras} (also known as {\it quiver Hecke algebras, or KLR algebras}) were introduced independently by Khovanov and Lauda \cite{KL2009,KL2011} and by Rouquier \cite{Rou}. In the Kac-Moody setting, the category of finitely generated graded projective modules over KLR algebras provides a categorification of the corresponding quantum groups. For symmetric Cartan data, the indecomposable projective modules correspond to Lusztig's canonical basis \cite{Rou1, VV2011}. Moreover, the cyclotomic quotients of KLR algebras categorify the irreducible highest weight representations of quantum groups \cite{KK2012}, as well as their crystals in the sense of Kashiwara \cite{LV2011}.

In the Kac-Moody case, Varagnolo and Vasserot \cite{VV2011} showed that KLR algebras admit a geometric realization as Steinberg-type convolution algebras associated with quiver flag varieties without loops, generalizing classical Springer theory. For quivers with arbitrary loops, a presentation of these Steinberg-type algebras was given in \cite{KKP2013}. Independently, Sauter \cite{Ju2013} developed a generalized quiver-graded Springer theory, in which the construction of \cite{KKP2013} appears as a special case for the type $A$ Weyl group.

In this setting, it is natural to expect that such generalized KLR algebras categorify quantum groups associated with quivers that may contain loops, and that they are closely related to Lusztig's theory of canonical bases. A partial link in this direction was established in \cite{KKP2013}, under restrictive assumptions on the quiver that exclude, in particular, the Jordan quiver. The quantum group considered there is relatively small, namely the {\it quantum generalized Kac-Moody algebra}.

The purpose of this paper is to develop an appropriate quantum group framework for these Steinberg-type algebras, thereby filling the gap left in \cite{KKP2013}.



T. Bozec studied in \cite{Bozec2014b} the theory of perverse sheaves on quiver representation varieties with possible loops, and proved that the Grothendieck group generated by Lusztig sheaves is spanned by elementary simple perverse sheaves, answering a question of Lusztig \cite{Lus93}.
In this framework, all partial flag varieties are taken into account, leading to a quantum group of Borcherds-Cartan type determined by a symmetric Borcherds-Cartan matrix $A=(a_{ij})_{i,j\in I}$, whose diagonal entries may be non-positive.
For an imaginary index $i$ with $a_{ii}\leq 0$, infinitely many generators $F_{in}$ ($n\in\Z_{>0}$) appear.

In the present paper, we restrict attention to the subalgebra $U^-$ of Bozec's quantum group arising from complete flag varieties $\FF_\nu$ for $\nu \in\N[I]$.
As a consequence, infinitely many generators $F_{i n}$ occur only in the isotropic case $a_{ii}=0$.

We first categorify $U^-$ by constructing the corresponding KLR algebras $R(\nu)$ for $\nu\in\N[I]$, using braid-like planar diagrammatics.
Rather than starting from geometric presentations of Steinberg algebras, which are technically more involved, we adopt an algebraic approach.
In particular, for $\nu =ni$ with $a_{ii}=0$, the algebra $R(ni)$ is the wreath product $\C[x_1,\ds,x_n]\rtimes \C[S_n]$, which coincides with the equivariant Borel--Moore homology ring $H_*^{GL_n}(\ZZ)$ of the Steinberg variety $\ZZ$ associated with the complete flag variety $GL_n/B_n$.

Let $R=\bigoplus_{\nu\in\N[I]} R(\nu)$, and let $K_0(R)$ denote the Grothendieck group of the category of finitely generated graded projective $R$-modules.
Using the framework developed in \cite{KL2009,KL2011}, together with a detailed analysis of the Jordan quiver case, we prove that there is a bialgebra isomorphism
$$\Ga: {_\A}U^- \xr{\sim} K_0(R),$$
where $\A=\Z[q,q^{-1}]$ and ${_\A}U^-$ denotes the $\A$-form of $U^-$.
As mentioned above, the Steinberg-type varieties $\ZZ_\nu$ associated with the complete flag varieties $\FF_\nu$ were studied in \cite{KKP2013,Ju2013}.
We show that $R(\nu)$ is isomorphic to the corresponding Steinberg-type convolution algebra $\R_\nu$.
As a consequence, combining classical Springer theory for $GL_n$ with the results of Varagnolo and Vasserot, we prove that the canonical basis elements of ${_\A}U^-$ correspond, under $\Ga$, to the indecomposable projective $R$-modules.

Finally, let $U$ be the quantum group obtained from $U^-$ via the Drinfeld double construction.
We show that the cyclotomic quotients $R^\Lambda$ for $\Lambda\in P^+$ provide a categorification of the irreducible highest weight $U$-module $V(\Lambda)$. The proof proceeds by iterating the lower-rank arguments developed in \cite{KK2012}.

\vskip 1mm
\noindent\textbf{Acknowledgements.}
  Young Rock Kim and Bolun Tong were supported by the National Research Foundation of Korea (NRF) grant funded by the Korea government (MSIT) (No. 2021R1A2C1011467). Young Rock Kim was supported by Hankuk University of Foreign Studies Research Fund.

\section{\textbf{Quantum groups of Borcherds-Cartan type}}

We introduce a class of quantum groups associated with quivers that may contain loops. These algebras appear as subalgebras of the quantum groups constructed in \cite{Bozec2014b} and admit a geometric realization in terms of complete flag varieties associated with arbitrary quivers.

\subsection{The algebra $U^-$}\

Let $I$ be a finite index set. An integer-valued symmetric matrix
$A=(a_{ij})_{i,j\in I}$ is called a {\it Borcherds-Cartan matrix} if it satisfies
\begin{itemize}
\item [(i)] $a_{ii}=2,0,-2,-4,\ds$,
\item [(ii)] $a_{ij}=a_{ji}\in \Z_{\leq 0}$ for $i\neq j$.
\end{itemize}

We define the following subsets of $I$:
$$I^+=\{i \in I \mid a_{ii}=2\},\quad
I^0=\{i \in I \mid a_{ii}=0\},\quad
I^-=\{i \in I \mid a_{ii}<0\},$$
and set $I^{\leq 0}= I^0 \sqcup I^-$. We further introduce the index set
$$\II=I^+ \sqcup I^- \sqcup \big( I^0 \ti \Z_{>0} \big).$$

Define a symmetric bilinear form
$\nu,\nu'\ma \nu\ac \nu'$
on the free abelian group $\Z[I]$ by
$$i\ac j=a_{ij}=a_{ji} \ \ \tx{for all}\ i,j\in I.$$
The triple $(I,A,\ac)$ is called a {\it Borcherds-Cartan datum}.

Let $q$ be an indeterminate. For $n,k\in \N$, set
$$[n]=\f{q^n-q^{-n}}{q-q^{-1}}, \quad [n]!=[n][n-1]\cdots [1], \quad {\begin{bmatrix} n \\ k \end{bmatrix}}=\f{[n]!}{[k]![n-k]!}.$$

\begin{definition}\label{U-} The quantum group $U^-$ associated to
 $(I, A, \ac)$  is the associative algebra over $\Q(q)$ with $1$ generated by $F_{i n}$ $((i,n)\in \II)$ satisfying the following relations:
$$\begin{aligned}
& \s_{r+s=1-n a_{ij}}(-1)^rF_i^{(r)}F_{jn}F_i^{(s)}=0 \quad \tx{for} \ i\in I^+,(j,n)\in \II \ \tx{and} \ i \neq (j,n),\\
& F_{in}F_{jm}-F_{jm}F_{in}=0 \quad\ \tx{for}\ a_{ij}=0.
\end{aligned}$$
Here we denote $F_i^{(n)}=F_i^n /[n]!$ for $i\in I^+$ and $n>0$. The algebra $U^-$ is $\N[I]$-graded by assigning $|F_{in}|=ni$.
\end{definition}

Define a twisted multiplication on $U^-\ot U^-$ by
$$(x_1\ot x_2)(y_1\ot y_2)=q^{-|x_2|\ac |y_1|}x_1y_1\ot x_2y_2,$$
for homogeneous  $x_1,x_2,y_1,y_2\in U^-$.
We have an algebra homomorphism $\rho: U^-\ra U^-\ot U^-$ (with respect to the twisted multiplication on $U^-\ot U^-$) given by
$$\rho(F_{in})=\s_{r+t=n} F_{ir}\ot F_{it} \ \ \tx{for}  \ (i,n)\in \II.$$
There is a nondegenerate symmetric bilinear form $\{ \ , \ \}: U^-\ti U^-\ra\Q(q)$ determined by
\begin{itemize}
\item[(1)] $\{x, y\}=0$ \ if $|x|\neq |y|$,
\item[(2)] $\{1,1\}=1$,
\item[(3)] $\{F_{in}, F_{in}\}=\ep_n:=\p_{k=1}^n\f{1}{1-q^{2k}}$ \ for  $(i,n)\in \II$,
\item[(4)] $\{x, yz\}=\{\rho(x), y\ot z\}$ \ for $x,y,z\in U^-$.
\end{itemize}

Let $\A=\Z[q,q^\zz]$ be the ring of Laurent polynomials. The $\A$-form $_{\A}U^-$ is the $\A$-subalgebra of $U^-$ generated by $F_i^{(n)}$  $(i\in I^+,n>0)$, $F_i$ $(i\in I^-)$ and  $F_{in}$ $(i\in I^0,n>0)$.

We denote by $^-$ the $\Q$-algebra involution of $U^-$ given by $\ovv{q}=q^\zz$ and
$\ovv{F_{in}}=F_{in}$.

\subsection{Geometric setting for $U^-$}\

\vs{1mm}

We briefly review the geometric construction for $U^-$ given in \cite{Bozec2014b, LL09, Lus93}. Let $(I,H)$ be a quiver, possibly with loops, where $I$ denote the set of vertices and $H$ the set of arrows. For each $h\in H$, let $h'$ and $h''$ denote its source and the target, respectively. Define
$$ h_i=\#\{\tx{loops on}\ i\}\ \ \tx{for}\ i\in I,\quad h_{ij}=\#\{h\in H\mid h'=i,h''=j\}\ \ \tx{for}\ i\neq j.$$

Fix $\nu=\sum_{i\in I}\nu_i i\in \N[I]$ with $\tx{ht}{(\nu)}:=\sum_{i\in I}\nu_i=n$. We set
$$V_\nu=\op_{i\in I}\C^{\nu_i},\quad
\EE_\nu=\op_{h\in H}\Hom(\C^{\nu_{h'}},\C^{\nu_{h''}}),\quad
G_{\nu}=\p_{i\in I}\tx{GL}_{\nu_i}(\C).$$
Then $G_{\nu}$ acts on $\EE_{\nu}$ by conjugation:
$$g\ac(x_h)=(g_{h''}x_h g_{h'}^\zz).$$

We denote by $D_{G_{\nu}}(\EE_{\nu})$ the bounded $G_\nu$-equivariant derived category of $\C$-constructible complexes on $\EE_{\nu}$, and by $P_{G_{\nu}}(\EE_{\nu})$ the abelian subcategory of $G_\nu$-equivariant perverse sheaves.

Let $\Seq(\nu)$ be the set of all sequences $\ii=i_1i_2\ds i_n$ in $I^n$ such that $\nu=i_1+i_2\cs+i_n$.
For $\ii \in \Seq(\nu)$, define
$$\begin{aligned}
& \FF_{\ii}=\{W =(0\subsetneq W_1\subsetneq \cs \subsetneq W_n=V_\nu)
\mid \underline{\tx{dim}}(W_k/W_{k-1})= i_k\}, \\
& \wi{\FF}_{\ii}=\{(x,W)\in \EE_{\nu}\ti\FF_{\ii}\mid x(W_k)\se W_{k-1}\}.
\end{aligned}$$

Let $G_\nu$ acts on $\wi{\FF}_{\ii}$ diagonally. The first projection $\pi_{\ii}:\wi{\FF}_{\ii}\ra \EE_{\nu}$ is a $G_{\nu}$-equivariant proper map, which yields $\LL_{\ii}=({\pi_{\ii}})_! (\C_{\wi{\FF}_{\ii}}[\tx{dim}\wi{\FF}_{\ii}])$  a semisimple complex in $ D_{G_{\nu}}(\EE_{\nu})$.

We denote:
\begin{itemize}
\item[(i)] $\PP_{\nu}$ is the set of isomorphism classes of simple perverse sheaves appearing  in $\LL_{\ii}$ as a summand with a possible shift for all $\ii\in\Seq(\nu)$,
\item[(ii)] $\QQ_{\nu}$ is the full subcategory of $D_{G_{\nu}}(\EE_{\nu})$ whose objects are finite direct sums of shifts of the simple perverse sheaves coming from $\PP_{\nu}$,
\item[(iii)] $\KK=\op_{\nu\in \N[I]}\KK_\nu$, where $\KK_\nu$ is the Grothendieck group of $\QQ_{\nu}$.
\end{itemize}

By \cite{Bozec2014b}, the group $\KK$ has a geometrically defined (twisted) $\A$-bialgebra structure that is isomorphic to the $_{\A}U^-$, associated with the Borcherds-Cartan matrix $A$ given by
\begin{equation} \label{quiver}
a_{ij} =
\begin{cases}
\ 2 - 2h_{i}  \ \ & \tx{for}\ i=j, \\
-h_{ij}-h_{ji} \ \ & \tx{for}\ i\neq j.
\end{cases}
\end{equation}
This isomorphism is given explicitly as follows:
\begin{equation*}
\begin{aligned}
& F_i^{(n)}  \longleftrightarrow {\C}_{\EE_{ni}} \ \ \tx{for}\ i \in I^+, \\
& F_i \longleftrightarrow  \C_{\EE_i} \ \ \tx{for}\ i \in I^-, \\
& F_{in} \longleftrightarrow  {\C}_{\{0\}\se \EE_{ni}} \ \ \tx{for}\  i\in I^0, \ n>0.
\end{aligned}
\end{equation*}
We identify $\KK$ with $_{\A}U^-$ via the isomorphism and refer to $\PP_{\nu}$ as the {\it canonical basis} of $_{\A}U^-_\nu$.

\begin{remark}\label{iso}
Let $i\in I^0$ and $\nu=ni$. We denote by $\{\OO_{\la}\}_{\la\vdash n}$ the nilpotent orbits (labelled by partitions of $n$) in $\EE_{ni}^{\tx{nil}}$ under the action of $GL_{ni}$. Then $F_{in}=IC(\OO_{(1^{n})})$ is the simple perverse sheaf associated with the closed orbit $\{0\}$  in $\EE_{ni}^{\tx{nil}}$.

Let $i\in I^0$. We abbreviate $ii\cs i\in\Seq(ni)$ by $i^n$ and write $F_i$ for $F_{(i,1)}$. The $n$-th power $F_i^n=\LL_{i^n}$ is the Springer sheaf
$$\Spr_{n}=\pi_!({\C}_{\wi{\NN}_n}[\tx{dim}\wi{\NN}_n]),$$
where $\NN_n=\EE_{ni}^{\tx{nil}}$, $\wi{\NN}_n=\wi{\FF}_{i^n}$, and $\pi=\pi_{i^n}:\wi{\NN}_n\ra \NN_n$ is the Springer map.

Therefore, we can write $F_i^n=\op_{\la\vdash n} IC(\OO_{\la})\ot V_{\la}$ for some nonzero vector spaces $V_{\la}$, and $F_i^n$ corresponding to the regular $\C[S_n]$-module under the Springer correspondence.
\end{remark}

There is a geometric pairing $\{ \ , \ \}:\KK\ti\KK\ra \Z(\!(q)\!)$
induced by the equivariant cohomology (see e.g. \cite[8.1.9]{Lus}), which coincide with the one we define on $_{\A}U^-$. In particular, for any $i\in I^0, n>0$, we have
$$\{F_{in},F_{in}\}=\s_j\tx{dim}\big(H^{j}_{GL_{n}}{(\tx{pt})}\big) q^j=\p_{k=1}^{n}\f{1}{1-q^{2k}}.$$

\section{\textbf{Categorification of ${U}^-$}}

\subsection{$\Z$-grading}\

\vs{1mm}

Let $R$ be a $\Z$-graded $\C$-algebra. For a graded $R$-module $M=\op_{n\in\Z}M_n$, its graded dimension is defined to be
$$\Dim M=\s_{n\in\Z}\big(\tx{dim}_{\C}M_n\big)q^n,$$
where $q$ is a formal variable. For $m\in \Z$, we denote by $M\{m\}$ the graded $R$-module obtained from $M$ by putting $(M\{m\})_n=M_{n-m}$. More generally, for $f(q)=\sum_{m\in \Z}a_mq^m\in\N[q,q^{-1}]$, we set $M^{f}=\op_{m\in\Z}(M\{m\})^{\oplus a_m}$.

Given two graded $R$-modules $M$ and $N$, we denote by $\Hom_{R\tx{-gr}}(M,N)$ the ${\C}$-vector space of degree-preserving homomorphisms and form the $\Z$-graded vector space
$$\HOM_R(M,N)=\op_{n\in\Z}\Hom_{R\tx{-gr}}(M\{n\},N)=\op_{n\in\Z}\Hom_{R\tx{-gr}}(M,N\{-n\}).$$

\subsection{KLR algebra $R(\nu)$}\

\vs{1mm}

Let $A$ be the Borcherds-Cartan matrix given by quiver as in \eqref{quiver}.
For $i\in I^{\leq 0}$, define the polynomial
$$H_i(u,v)=(-1)^{a_{ii}/2}{(u-v)}^{-a_{ii}}.$$
For $i,j\in I$, define
$$Q_{ij}(u,v)=(-1)^{h_{ij}}{(u-v)}^{-a_{ij}}.$$

Fix $\nu\in\N[I]$ with $\tx{ht}{(\nu)}=n$. The KLR algebra
$R(\nu)$ associated to $A$ is defined to be the $\C$-algebra of the homogeneous generators given by diagrams (see \cite{KL2009} for a detailed explanation of the braid-like planar diagrams):

$$\begin{aligned}
& 1_{\ii} =\ {\fontsize{11}{11}\selectfont\xy
(0,6)*{}; (0,-6)*{} **\dir{-};
(5,0)*{\cs};
(10,6)*{}; (10,-6)*{} **\dir{-};
(15,0)*{\cs};
(20,6)*{}; (20,-6)*{} **\dir{-};
(0.7,-8.7)*{i_1}; (10.7,-8.7)*{i_k}; (20.7,-8.7)*{i_n};
\endxy} \quad\ \tx{for} \ \ii \in \Seq(\nu), \ \ \tx{deg}(1_{\ii})=0, \\ \\
& x_{k,\ii} =\ {\fontsize{11}{11}\selectfont\xy
(0,6)*{}; (0,-6)*{} **\dir{-};
(5,0)*{\cs};
(10,6)*{}; (10,-6)*{} **\dir{-};
(10,0)*{\bu}; (15,0)*{\cs};
(20,6)*{}; (20,-6)*{} **\dir{-};
(0.7,-8.7)*{i_1}; (10.7,-8.7)*{i_k}; (20.7,-8.7)*{i_n};
\endxy} \quad\ \tx{for} \ \ii\in \Seq(\nu),\ 1\leq k\leq n, \ \ \tx{deg}(x_{k,\ii})=2,\\ \\
& \tau_{k,\ii}=\ {\fontsize{11}{11}\selectfont\xy
(0,6)*{}; (0,-6)*{} **\dir{-};
(4,0)*{\cs};
(6.5,6)*{}; (14.5,-6)*{} **\dir{-};
(14.5,6)*{}; (6.5,-6)*{} **\dir{-};
(17,0)*{\cs};
(21,6)*{}; (21,-6)*{} **\dir{-};
(0.7,-8.7)*{i_1}; (7,-8.7)*{i_k}; (14.7,-8.7)*{i_{k+1}}; (21.7,-8.7)*{i_n};
\endxy}  \quad\ \tx{for} \ \ii\in \Seq(\nu),\ 1\leq k\leq n-1, \ \ \tx{deg}(\tau_{k,\ii})=-i_k\ac i_{k+1}.
\end{aligned}$$

\vs{3mm}

\noindent Subject to the following local relations:

\begin{equation}
{\fontsize{10}{10}\selectfont\xy
(0,7)*{}; (0,-7)*{} **\crv{(11,0)};
(7,-7)*{}; (7,7)*{} **\crv{(-4,0)};
(0,-9)*{i}; (7,-9)*{j};
\endxy} \  =  \
\begin{cases}
 \ \ 0 \ & \tx{ if}\ i= j\in I^+, \\  \\
 \   H_i\Big(\ {\fontsize{10}{10}\selectfont\xy
(0,4)*{}; (0,-4)*{} **\dir{-};
(5,4)*{}; (5,-4)*{} **\dir{-};
(0,0)*{\bu}; (0,-6)*{i}; (5,-6)*{i};
\endxy} \ , \ {\fontsize{10}{10}\selectfont\xy
(0,4)*{}; (0,-4)*{} **\dir{-};
(5,4)*{}; (5,-4)*{} **\dir{-};
(5,0)*{\bu}; (0,-6)*{i}; (5,-6)*{i};
\endxy}\ \Big)  & \tx{ if}\ i= j\in I^{\leq 0}, \\ \\
\ Q_{ij}\Big(\ {\fontsize{10}{10}\selectfont\xy
(0,4)*{}; (0,-4)*{} **\dir{-};
(5,4)*{}; (5,-4)*{} **\dir{-};
(0,0)*{\bu}; (0,-6)*{i}; (5,-6)*{j};
\endxy} \ , \ {\fontsize{10}{10}\selectfont\xy
(0,4)*{}; (0,-4)*{} **\dir{-};
(5,4)*{}; (5,-4)*{} **\dir{-};
(5,0)*{\bu}; (0,-6)*{i}; (5,-6)*{j};
\endxy}\ \Big) & \tx{ if}\ i\neq  j,
\end{cases}
\end{equation}

\begin{equation}
{\fontsize{10}{10}\selectfont\xy
(0,5)*{}; (8,-5)*{} **\dir{-}?(.25)*{\bu};
(8,5)*{}; (0,-5)*{} **\dir{-};
(0,-7)*{i}; (8,-7)*{i}; \endxy} \ -\
{\fontsize{10}{10}\selectfont\xy
(0,5)*{}; (8,-5)*{} **\dir{-}?(.75)*{\bu};
(8,5)*{}; (0,-5)*{} **\dir{-};
(0,-7)*{i}; (8,-7)*{i}; \endxy}\ =\
{\fontsize{10}{10}\selectfont\xy
(0,5)*{}; (0,-5)*{} **\dir{-};
(6,5)*{}; (6,-5)*{} **\dir{-};
(0,-7)*{i}; (6,-7)*{i};
\endxy}
\quad\quad\quad
{\fontsize{10}{10}\selectfont\xy
(0,5)*{}; (8,-5)*{} **\dir{-};
(8,5)*{}; (0,-5)*{} **\dir{-}?(.75)*{\bu};
(0,-7)*{i}; (8,-7)*{i}; \endxy}  \ - \
{\fontsize{10}{10}\selectfont\xy
(0,5)*{}; (8,-5)*{} **\dir{-};
(8,5)*{}; (0,-5)*{} **\dir{-}?(.25)*{\bu};
(0,-7)*{i}; (8,-7)*{i}; \endxy}\ =\
{\fontsize{10}{10}\selectfont\xy
(0,5)*{}; (0,-5)*{} **\dir{-};
(6,5)*{}; (6,-5)*{} **\dir{-};
(0,-7)*{i}; (6,-7)*{i};
\endxy}  \quad\ \tx{if}\ i \in I^+,
\end{equation}

\begin{equation}
{\fontsize{10}{10}\selectfont\xy
(0,5)*{}; (8,-5)*{} **\dir{-}?(.25)*{\bu};
(8,5)*{}; (0,-5)*{} **\dir{-};
(0,-7)*{i}; (8,-7)*{j}; \endxy} \   =  \ {\fontsize{10}{10}\selectfont\xy
(0,5)*{}; (8,-5)*{} **\dir{-}?(.75)*{\bu};
(8,5)*{}; (0,-5)*{} **\dir{-};
(0,-7)*{i}; (8,-7)*{j}; \endxy}
\quad\quad\quad
{\fontsize{10}{10}\selectfont\xy
(0,5)*{}; (8,-5)*{} **\dir{-};
(8,5)*{}; (0,-5)*{} **\dir{-}?(.75)*{\bu};
(0,-7)*{i}; (8,-7)*{j}; \endxy}  \ =  \ {\fontsize{10}{10}\selectfont\xy
(0,5)*{}; (8,-5)*{} **\dir{-};
(8,5)*{}; (0,-5)*{} **\dir{-}?(.25)*{\bu};
(0,-7)*{i}; (8,-7)*{j}; \endxy}  \quad\ \tx{otherwise},
\end{equation}

\begin{equation}
{\fontsize{10}{10}\selectfont\xy
(0,7)*{}; (12,-7)*{} **\crv{(0,0)&(12,0)}; (12,7)*{}; (0,-7)*{} **\crv{(12,0)&(0,0)};
(6,7)*{}; (6,-7)*{} **\crv{(6,6)&(3,4.5)&(0,1.5)&(0,-1.5)&(3,-4.5)&(6,-6)};
(0,-9)*{i}; (6,-9)*{j};(12,-9)*{i};
\endxy} \ - \ {\fontsize{10}{10}\selectfont\xy
(0,7)*{}; (12,-7)*{} **\crv{(0,0)&(12,0)}; (12,7)*{}; (0,-7)*{} **\crv{(12,0)&(0,0)};
(6,7)*{}; (6,-7)*{} **\crv{(6,6)&(8,4.5)&(12,1.5)&(12,-1.5)&(8,-4.5)&(6,-6)};
(0,-9)*{i}; (6,-9)*{j};(12,-9)*{i};
\endxy}  \ =\ \frac{Q_{ij}\Big(\ { \fontsize{8}{8}\selectfont\xy
(0,3.5)*{}; (0,-3.5)*{} **\dir{-};
(3,3.5)*{}; (3,-3.5)*{} **\dir{-};
(6,3.5)*{}; (6,-3.5)*{} **\dir{-};
(0,0)*{\bu}; (0,-5)*{i}; (3,-5)*{j};(6,-5)*{i};
\endxy}\ , \ {\fontsize{8}{8}\selectfont\xy
(0,3.5)*{}; (0,-3.5)*{} **\dir{-};
(3,3.5)*{}; (3,-3.5)*{} **\dir{-};
(6,3.5)*{}; (6,-3.5)*{} **\dir{-};
(3,0)*{\bu}; (0,-5)*{i}; (3,-5)*{j};(6,-5)*{i};
\endxy}\ \Big)-Q_{ij}\Big(\ { \fontsize{8}{8}\selectfont\xy
(0,3.5)*{}; (0,-3.5)*{} **\dir{-};
(3,3.5)*{}; (3,-3.5)*{} **\dir{-};
(6,3.5)*{}; (6,-3.5)*{} **\dir{-};
(6,0)*{\bu}; (0,-5)*{i}; (3,-5)*{j};(6,-5)*{i};
\endxy}
\ , \ {\fontsize{8}{8}\selectfont\xy
(0,3.5)*{}; (0,-3.5)*{} **\dir{-};
(3,3.5)*{}; (3,-3.5)*{} **\dir{-};
(6,3.5)*{}; (6,-3.5)*{} **\dir{-};
(3,0)*{\bu}; (0,-5)*{i}; (3,-5)*{j};(6,-5)*{i};
\endxy}\ \Big)}  { {\fontsize{8}{8}\selectfont\xy
(0,3.5)*{}; (0,-3.5)*{} **\dir{-};
(3,3.5)*{}; (3,-3.5)*{} **\dir{-};
(6,3.5)*{}; (6,-3.5)*{} **\dir{-};
(0,0)*{\bu}; (0,-5)*{i}; (3,-5)*{j};(6,-5)*{i};
\endxy}-{\fontsize{8}{8}\selectfont\xy
(0,3.5)*{}; (0,-3.5)*{} **\dir{-};
(3,3.5)*{}; (3,-3.5)*{} **\dir{-};
(6,3.5)*{}; (6,-3.5)*{} **\dir{-};
(6,0)*{\bu}; (0,-5)*{i}; (3,-5)*{j};(6,-5)*{i};
\endxy}}
\quad \tx{if}\ i\in I^+, i\ne j \ \tx{and} \ a_{ij} < 0,
\end{equation}

\begin{equation}
{\fontsize{10}{10}\selectfont\xy
(0,7)*{}; (12,-7)*{} **\crv{(0,0)&(12,0)}; (12,7)*{}; (0,-7)*{} **\crv{(12,0)&(0,0)};
(6,7)*{}; (6,-7)*{} **\crv{(6,6)&(3,4.5)&(0,1.5)&(0,-1.5)&(3,-4.5)&(6,-6)};
(0,-9)*{i}; (6,-9)*{j};(12,-9)*{k};
\endxy} \  = \ {\fontsize{10}{10}\selectfont\xy
(0,7)*{}; (12,-7)*{} **\crv{(0,0)&(12,0)}; (12,7)*{}; (0,-7)*{} **\crv{(12,0)&(0,0)};
(6,7)*{}; (6,-7)*{} **\crv{(6,6)&(8,4.5)&(12,1.5)&(12,-1.5)&(8,-4.5)&(6,-6)};
(0,-9)*{i}; (6,-9)*{j};(12,-9)*{k};
\endxy} \quad\quad \tx{otherwise}.
\end{equation}

\vs{3mm}

For $\ii,\jj\in \nu$, set $_{\jj}R(\nu)_{\ii}=1_{\jj}R(\nu)1_{\ii}$, then we have the decomposition
$$R(\nu)=\op_{\ii,\jj} {_{\jj}R(\nu)_{\ii}}.$$
We denote by $P_{\ii}=R(\nu)1_{\ii}$ (resp. $_{\jj}P=1_{\jj}R(\nu)$) the gr-projective left  $R(\nu)$-module (resp. right $R(\nu)$-module).

For $\ii \in \Seq(\nu)$, let $\PPP_{\ii}={\C}[x_1(\ii),\ds,x_n(\ii)]$ and form the $\C$-vector space $\PPP_{\nu}=\op_{\ii\in \Seq(\nu)}\PPP_{\ii}$. The symmetric group $S_n$ acts on $\PPP_{\nu}$ by sending $x_a(\ii)$ to $x_{\om(a)}(\om(\ii))$ for each $\om\in S_n$.

We define an action of $R({\nu})$ on $\PPP_{\nu}$ as follows:
\begin{itemize}
\item[(1)]  If $\ii\neq\kk$, then $_{\jj}R(\nu)_{\ii}$ acts on $\PPP_{\kk}$ by $0$.

\vs{1mm}

\item[(2)]  For $f\in \PPP_{\ii}$,
$1_{\ii}\ac f=f, \ x_{k,\ii}\ac f=x_k{(\ii)}f.$

\vs{1mm}

\item[(3)] If $\ii=i_1\ds i_n$ with $i_k=i$ and $i_{k+1}=j$, then
\begin{equation*}
\tau_{k,\ii}\ac f=\begin{cases}\pl_k(f)=\f{f- s_kf}{x_k{(\ii)}-x_{k+1}(\ii)} & \tx{if}\ i= j\in I^+,\\ \vs{-4.5mm} \\
{\left(x_k(\ii)-x_{k+1}(\ii)\right)}^{-a_{ii}/2}s_kf & \tx{if}\ i=j\in I^{\leq 0}, \\ \vs{-5.3mm} \\
{\left(x_k(s_k\ii)-x_{k+1}(s_k\ii)\right)}^{h_{ij}} s_kf & \tx{if}\ i\neq j.
\end{cases}
\end{equation*}
\end{itemize}
It is easy to check $\PPP_\nu$ is a well-defined $R(\nu)$-module.

\subsection{Algebras $R(ni)$ and their gr-irreducible modules}\

\vs{1mm}

Let $\nu=ni$ for some $i\in I$. Note that the unique sequence in $\Seq(ni)$ is $i^n$.
The  algebra $R(ni)$ is generated by $x_{1,i^n},\dots, x_{n,i^n}$ of degree $2$
and $\tau_{1,i^n},\dots,\tau_{n-1,i^n}$ of degree $-a_{ii}$
subject to the following local relations:

\vs{-3mm}

$${\fontsize{9}{9}\selectfont\xy
(0,6)*{}; (0,-6)*{} **\crv{(10, 0)};
(6,6)*{}; (6,-6)*{} **\crv{(-4,0)};
(0,-8)*{i}; (6,-8)*{i};
\endxy}=0 \quad\quad \
{\fontsize{9}{9}\selectfont\xy
(0,5)*{}; (8,-5)*{} **\dir{-}?(.25)*{\bu};
(8,5)*{}; (0,-5)*{} **\dir{-};
(0,-7)*{i}; (8,-7)*{i}; \endxy} -
{\fontsize{9}{9}\selectfont\xy
(0,5)*{}; (8,-5)*{} **\dir{-}?(.75)*{\bu};
(8,5)*{}; (0,-5)*{} **\dir{-};
(0,-7)*{i}; (8,-7)*{i}; \endxy}\ =\
{\fontsize{9}{9}\selectfont\xy
(0,5)*{}; (8,-5)*{} **\dir{-};
(8,5)*{}; (0,-5)*{} **\dir{-}?(.75)*{\bu};
(0,-7)*{i}; (8,-7)*{i}; \endxy} -
{\fontsize{9}{9}\selectfont\xy
(0,5)*{}; (8,-5)*{} **\dir{-};
(8,5)*{}; (0,-5)*{} **\dir{-}?(.25)*{\bu};
(0,-7)*{i}; (8,-7)*{i}; \endxy}\ =\
{\fontsize{9}{9}\selectfont\xy
(0,5)*{}; (0,-5)*{} **\dir{-};
(6,5)*{}; (6,-5)*{} **\dir{-};
(0,-7)*{i}; (6,-7)*{i};
\endxy} \quad\quad \
{\fontsize{9}{9}\selectfont\xy
(0,6)*{}; (10,-6)*{} **\crv{(0,0)&(10,0)}; (10,6)*{}; (0,-6)*{} **\crv{(10,0)&(0,0)};
(5,6)*{}; (5,-6)*{} **\crv{(5,5)&(3,4)&(0,1.5)&(0,-1.5)&(3,-4)&(5,-5)};
(0,-8)*{i}; (5,-8)*{i};(10,-8)*{i};
\endxy} =
{\fontsize{9}{9}\selectfont\xy
(0,6)*{}; (10,-6)*{} **\crv{(0,0)&(10,0)}; (10,6)*{}; (0,-6)*{} **\crv{(10,0)&(0,0)};
(5,6)*{}; (5,-6)*{} **\crv{(5,5)&(7,4)&(10,1.5)&(10,-1.5)&(7,-4)&(5,-5)};
(0,-8)*{i}; (5,-8)*{i};(10,-8)*{i};
\endxy} \quad\quad \ \tx{if}\ i \in I^+ .$$
$${\fontsize{9}{9}\selectfont\xy
(0,6)*{}; (0,-6)*{} **\crv{(10, 0)};
(6,6)*{}; (6,-6)*{} **\crv{(-4,0)};
(0,-8)*{i}; (6,-8)*{i};
\endxy}=(-1)^{\frac{a_{ii}}{2}}\Big(\ { \fontsize{9}{9}\selectfont\xy
(0,5)*{}; (0,-5)*{} **\dir{-};
(5,5)*{}; (5,-5)*{} **\dir{-};
(0,0)*{\bu}; (0,-7)*{i}; (5,-7)*{i};
\endxy}\ -\
{\fontsize{9}{9}\selectfont\xy
(0,5)*{}; (0,-5)*{} **\dir{-};
(5,5)*{}; (5,-5)*{} **\dir{-};
(5,0)*{\bu};  (0,-7)*{i}; (5,-7)*{i};
\endxy} \ \ \Big)^{-a_{ii}}\quad\
{\fontsize{9}{9}\selectfont\xy
(0,5)*{}; (8,-5)*{} **\dir{-}?(.25)*{\bu};
(8,5)*{}; (0,-5)*{} **\dir{-};
(0,-7)*{i}; (8,-7)*{i}; \endxy} =
{\fontsize{10}{10}\selectfont\xy
(0,5)*{}; (8,-5)*{} **\dir{-}?(.75)*{\bu};
(8,5)*{}; (0,-5)*{} **\dir{-};
(0,-7)*{i}; (8,-7)*{i}; \endxy}  \quad
{\fontsize{10}{10}\selectfont\xy
(0,5)*{}; (8,-5)*{} **\dir{-};
(8,5)*{}; (0,-5)*{} **\dir{-}?(.75)*{\bu};
(0,-7)*{i}; (8,-7)*{i}; \endxy}   =
{\fontsize{10}{10}\selectfont\xy
(0,5)*{}; (8,-5)*{} **\dir{-};
(8,5)*{}; (0,-5)*{} **\dir{-}?(.25)*{\bu};
(0,-7)*{i}; (8,-7)*{i}; \endxy} \quad \ \ \
{\fontsize{9}{9}\selectfont\xy
(0,6)*{}; (10,-6)*{} **\crv{(0,0)&(10,0)}; (10,6)*{}; (0,-6)*{} **\crv{(10,0)&(0,0)};
(5,6)*{}; (5,-6)*{} **\crv{(5,5)&(3,4)&(0,1.5)&(0,-1.5)&(3,-4)&(5,-5)};
(0,-8)*{i}; (5,-8)*{i};(10,-8)*{i};
\endxy} =
{\fontsize{9}{9}\selectfont\xy
(0,6)*{}; (10,-6)*{} **\crv{(0,0)&(10,0)}; (10,6)*{}; (0,-6)*{} **\crv{(10,0)&(0,0)};
(5,6)*{}; (5,-6)*{} **\crv{(5,5)&(7,4)&(10,1.5)&(10,-1.5)&(7,-4)&(5,-5)};
(0,-8)*{i}; (5,-8)*{i};(10,-8)*{i};
\endxy} \quad \tx{if}\ i \in I^{\leq 0}.$$

\vs{1mm}

We will abbreviate $x_{k,i^n}$ (resp. $\tau_{k,i^n}$) for $x_k$ (resp. $\tau_k$). In all cases, $R(ni)$ has a basis
$$\{x_1^{r_1}\cs x_n^{r_n}\tau_\om\mid\om\in S_n,\ r_1,\ds,r_n\geq 0\}.$$
Indeed, for example, if $i\in I^{\leq 0}$, it suffices to consider the actions of these elements on $x_1^Nx_2^N\cs x_n^{nN}$ for $N\gg 0$. Thus, we can identify the polynomial algebra $P_n=\C[x_1,\ds,x_n]$ with the subalgebra of $R(ni)$ generated by $x_1,\ds,x_n$. Consequently, the center of $R(ni)$ is identified with $Z_n$, the algebra of symmetric polynomials in $x_1,\ds,x_n$.

We next consider the graded irreducible representations of $R(ni)$.

{\bf Case 1}: $i\in I^+$.

By the representation theory of the nil-Hecke algebras, $R(ni)$ has a unique gr-irreducible module $V(i^n)$ of graded dimension $[n]!$, which is isomorphic to
$$R(ni)\ot_{P_n}\mathbf 1_n\left\{\f{n(n-1)}{2}\right\}.$$
Here $\mathbf 1_n$ is the one-dimensional trivial module over $P_n$ on which  each $x_k$ acts by $0$.

{\bf Case 2}: $i\in I^-$.

In this case, $R(ni)$ has only trivial idempotents. It has a unique gr-irreducible module $V(i^n)$, which is the one-dimensional trivial module with the gr-projective cover $R(ni)$.

{\bf Case 3}: $i\in I^0$.

In the case, the degree zero part $R(ni)_0$  is the symmetric group algebra $\C[S_n]$. Let $\Pp_n$ be the set of partitions of $n$. It is well known that $\C[S_n]$ has $|\Pp_n|$ many irreducible modules.

Let $V$ be an irreducible $\C[S_n]$-module (which is also an indecomposable projective module). Then, we define
$$\wi{V}:=\f{R(ni)\ot_{R(ni)_0}V} {R(ni)_{>0}\ot_{R(ni)_0}V} $$
as a gr-irreducible $R(ni)$-module. In other words, $\wi{V}$ is obtained from $V$ by annihilating the actions of $x_1,\ds,x_n$. Moreover, all gr-irreducible $R(ni)$-modules can be obtained in this way.

When no confusion arises, we denote $\wi{V}$ simply by $V$. We have shown the following:

\begin{proposition}
If $V_1,\ds,V_{|\Pp_n|}$ is a complete set of non-isomorphic classes of irreducible $\C S_n$-modules, then $V_1,\ds,V_{|\Pp_n|}$ is a complete set of non-isomorphic classes of gr-irreducible $R(ni)$-modules. In particular, the gr-Jacobson radical $J^{\tx{gr}}(R(ni))=R(ni)_{>0}$.
\end{proposition}

Let $V_{i,n}$ be the one-dimensional trivial module over $\C[S_n]$. Note that
$$V_{i,n}=\C[S_n]\ac e_{i,n}=\C\ac e_{i,n},$$
where $e_{i,n}$ is defined by
$$e_{i,n}=\f{1}{n!}\s_{\om \in S_n} \om\in\C[S_n].$$
If $r+t=n$, then the restriction of $V_{i,n}$ to $\C [S_r]\ot\C[S_t]$-modules gives
$$\Res^n_{r,t}V_{i,n}\co V_{i,r}\ot V_{i,t}=\C[S_r]\ac e_{i,r}\ot \C[S_t]\ac e_{i,t}.$$

The gr-projective cover of $R(ni)$-module $V_{i,n}$ is $P_{i,n}=R(ni)e_{i,n}$, which has a basis
$$\{x_1^{r_1}\cs x_n^{r_n}\ac e_{i,n}\mid r_1,\ds,r_n\geq 0\}.$$
The restriction of $P_{i,n}$ to $R(ri)\ot R(ti)$-modules gives
\begin{equation}\label{P}
\Res^n_{r,t}P_{i,n}\co P_{i,r}\ot P_{i,t}.
\end{equation}

Since $e_{i,n}R(ni)e_{i,n}$ is spanned by $\{f\ac e_{i,n}\mid f\in Z_n\}$, we see that
\begin{equation}\label{P1}
(P_{i,n},P_{i,n})= \Dim e_{i,n}R(ni)e_{i,n}= \Dim Z_n=\p_{k=1}^{n}\f{1}{1-q^{2k}}.
\end{equation}
Here $( \ , \ )$ is the Khovanov-Lauda's form defined in (\ref{KL-form}).

\subsection{Grothendieck groups $K_0(R)$ and $G_0(R)$}\

\vs{1mm}

Using the polynomial representation $\PPP_{\nu}$ of $R(\nu)$, one can obtain by a similar argument in \cite[Theorem 2.5]{KL2009} that $\PPP_{\nu}$ is a faithful $R(\nu)$-module and  $_{\ii}R(\nu)_{\jj}$ $(\ii,\jj\in \nu)$  has a basis
\begin{equation}\label{PBW}
\{x_{1,\ii}^{r_1}\cs x_{n,\ii}^{r_n}\ac\widehat\om_\jj \mid r_1,\ds,r_n\geq 0,\ \om\in S_n\ \tx{such that} \ \om(\jj)=\ii\},
\end{equation}
where $\widehat\om_\jj\in {_{\ii}} R(\nu)_{\jj}$ is uniquely determined by a fixed reduced expression of $\om$.

Assume that there is a sequence ${i_1}^{m_1}\cs{i_t}^{m_t}\in\Seq(\nu)$ such that $i_1,\ds,i_t$ are all distinct. Similar to \cite[Theorem 2.9]{KL2009}, the center $Z(R{(\nu)})$ of $R{(\nu)}$ can be described as
$$Z(R(\nu))\co\bigotimes^t_{k=1}\C[z_{1},\ds,z_{m_k}]^{S_{m_k}},$$
where the latter is a tensor product of symmetric polynomial algebras such that the generators in $\C[z_1,\ds,z_{m_k}]$ are of degree $2$. Moreover, $R(\nu)$ is a free $Z(R(\nu))$-module of rank $\left((m_1+\cs+m_t)!\right)^2$. It is also a gr-free $Z(R(\nu))$-module of finite rank. So we have
$$\Dim Z(R(\nu))=\p^t_{k=1}\left(\p^{m_k}_{s=1}\f{1}{1-q ^{2s}}\right)$$
and $\Dim R(\nu)\in \Z[q,q^\zz]\ac\Dim Z(R(\nu))$.

Denote by
$$\begin{aligned}
&  R(\nu)\Mod: \tx{the category of finitely generated graded}\  R(\nu) \tx{-modules},\\
& R(\nu)\fMod: \tx{the category of finite-dimensional graded}\  R(\nu) \tx{-modules},\\
& R(\nu)\pMod: \tx{the category of projective objects in } R(\nu)\Mod.
\end{aligned}$$

Up to isomorphism and degree shifts, each $R(\nu)$ has only finitely many gr-irreducible modules, all of which are finite-dimensional and are irreducible $R(\nu)$-modules by forgetting the grading.

Let $\B_{\nu}$ be the set of equivalence classes of gr-irreducible $R(\nu)$-modules. Choose one representative $S_b$ from each equivalence class and denote by $P_b$ the gr-projective cover of $S_b$. The Grothendieck group $G_0(R(\nu))$ (resp. $K_0(R(\nu))$) of $R(\nu)\fMod$ (resp. $R(\nu)\pMod$) are free $\Z[q,q^\zz]$-modules with $q[M]=[M\{1\}]$, and with a basis $\{[S_b]\}_{b\in \B_{\nu}}$ (resp. $\{[P_b]\}_{b\in \B_{\nu}}$).

Let $R=\op_{\nu\in \N[I]}R(\nu)$ and form
$$G_0(R)=\op_{\nu\in\N[I]}G_0(R(\nu)),\quad K_0(R)=\op_{\nu\in\N[I]}K_0(R(\nu)).$$

The $\Z[q,q^\zz]$-modules $K_0(R)$ and $G_0(R)$ are equipped with twisted bialgebras structure induced by the induction and restriction functors
\begin{equation*}
\begin{aligned}
&\Ind^{\nu+\nu'}_{\nu,\nu'}: R(\nu)\ot R(\nu')\Mod\ra R(\nu+\nu')\Mod,\ M\ma R(\nu+\nu')1_{\nu,\nu'}\ot_{R(\nu)\ot R(\nu')}M,\\
&\Res^{\nu+\nu'}_{\nu,\nu'}: R(\nu+\nu')\Mod\ra R(\nu)\ot R(\nu')\Mod,\ N\ma 1_{\nu,\nu'}N,
\end{aligned}
\end{equation*}
where $1_{\nu,\nu'}=1_\nu\ot 1_{\nu'}$.

More precisely, we set $|x|=|\nu|\in \N[I]$ for $x\in R(\nu)\Mod$ and equip $K_0(R)\ot K_0(R)$ (resp. $G_0(R)\ot G_0(R)$) with a twisted algebra structure via
$$(x_1\ot x_2)(y_1\ot y_2)=q^{-|x_2|\ac |y_1|}x_1y_1\ot x_2y_2,$$
then $\Res$ is a $\Z[q,q^\zz]$-algebra homomorphism by Mackey's Theorem \cite[Proposition 2.18]{KL2009}.

The $K_0(R)$ and $G_0(R)$ are dual to each other with respect to the bilinear pairing
$( \ , \ ): K_0(R)\ti G_0(R)\ra \Z[q,q^\zz]$ defined by
\begin{equation}\label{KL-form}([P],[M])=\Dim P^\psi\ot_{R(\nu)}M=\Dim \HOM_{R(\nu)}(\ovv{P},M),\end{equation}
where $\psi$ is the anti-involution of $R(\nu)$ obtained by flipping diagrams about horizontal axis, which turns a left $R(\nu)$-module into a right one, and $\ovv{P}=\HOM(P,R(\nu))^\psi$.

There is also a symmetric bilinear form $( \ , \ ): K_0(R)\ti K_0(R)\ra \Z(\!(q)\!)$ defined in the same way. By (\ref{P1}) and \cite[Proposition 3.3]{KL2009}, we obtain the following proposition.

\begin{proposition}\label{Sy} The symmetric bilinear form on $K_0(R)$ satisfies:
\begin{itemize}
\item[(1)] $([M],[N])=0$ \ if $M\in R(\nu)\Mod, N\in R(\mu)\Mod$ with $\nu\neq\mu$.
\item[(2)] $(1,1)=1$, where $1=\C$ as a module over $R(0)=\C$.
\item[(3)] $([P_{i}], [P_{i}]) =1/(1-q^2)$ \ \ for $i\in I^\pm$,\\
$([P_{i,n}], [P_{i,n}]) =1/\left((1-q^2)(1-q^4)\cs (1-q^{2n})\right)$ \ \ for $i\in I^0,n>0$.
\item[(4)] $(x, yz) = (\Res(x), y \otimes z)$  \ for $x,y,z\in K_0(R)$.
\end{itemize}
\end{proposition}

\subsection{Quantum Serre relations}\

\vs{1mm}

Let $\ii$ be a sequence with {\it divided powers} in $\nu$:
$$\ii=j_1\ds j_{p_0}i_1^{(m_1)}k_1\ds k_{p_1}i_2^{(m_2)}\ \ds\ i_t^{(m_t)} h_1\ds h_{p_t}, $$
where $i_1,\ds,i_t\in I^+$ and $\ii$ is of the weight $\nu$.

For such an $\ii$, we assign the following idempotent of $R(\nu)$
$$1_{\ii}=1_{j_1\ds j_{p_0}}\ot e_{i_1,m_1}\ot 1_{k_1\ds k_{p_1}}\ot e_{i_2,m_2} \ot\ \cs\ \ot e_{i_t,m_t}\ot 1_{h_1\ds h_{p_t}},$$
where $e_{i,m}=x_1^{m-1}x_2^{m-2}\cs x_{m-1}\tau_{w_0}$ with $w_0$ being the longest element in $S_{m}$.

We define the gr-projective modules for $\ii$ by
$$ {_\ii P} =1_\ii R(\nu)\{-\langle\ii \rangle\},\quad  \ P_\ii= R(\nu)\psi(1_\ii)\{-\langle\ii \rangle\},$$
where
$$\langle\ii\rangle=\s_{k=1}^t\f{m_k(m_k-1)}{2}.$$
In particular, for $i\in I^+$ and $n\geq 0$,
\begin{equation}\label{e}
P_{i^{(n)}}= R(ni)\psi{(e_{i,n})}\left\{-\f{n(n-1)}{2} \right\} \co R(ni)e_{i,n}\left\{ \f{n(n-1)}{2} \right\} .
\end{equation}

\begin{proposition}\label{Serre}
Suppose $i\in I^+, j\in I$ and $i\neq j$. Let $n\in \Z_{>0}$ and $m=1-na_{ij}$.
We have isomorphisms of graded left $R(\nu)$-modules
$$\begin{aligned}
&\op^{\lfloor \frac{m}{2} \rfloor}_{c=0}P_{i^{(2c)}j^ni^{(m-2c)}}\co\op^{\lfloor \frac{m-1}{2} \rfloor}_{c=0}P_{i^{(2c+1)}j^ni^{(m-2c-1)}}\quad\tx{if}\ j\in I^\pm,\\
& \op^{\lfloor \frac{m}{2} \rfloor}_{c=0}R(i^mj^n)\psi(e_{i,2c}\ot e_{j,n}\ot e_{i,m-2c})\{-\langle i^{(2c)}i^{(m-2c)} \rangle\}\\
&\phantom{aaa} \co\op^{\lfloor \frac{m-1}{2} \rfloor}_{c=0}R(i^mj^n)\psi(e_{i,2c+1}\ot e_{j,n}\ot e_{i,m-2c-1})\{-\langle i^{(2c+1)}i^{(m-2c-1)}\rangle\}\quad\tx{if}\ j\in I^0.
\end{aligned}$$
Moreover, if $a_{ij}=0$, then
$$\begin{aligned}
& P_{ij}\co P_{ji}\quad \tx{for}\ i,j\in I,\\
& R(i^nj)\ac e_{i,n}\ot 1_j\co  R(i^n j)\ac 1_j\ot e_{i,n}\quad \tx{for}\ i\in I^0, j\in I,\\
& R(i^n j^m)\ac e_{i,n}\ot e_{j,m}\co R(i^n j^m)\ac e_{j,m}\ot e_{i,n} \quad \tx{for}\ i,j \in I^0.
\end{aligned}$$
\begin{proof}
The proof is the same as the `Box' calculations in \cite{KL2011}. We only explain the last isomorphism. Let $i,j\in I^0$ with $a_{ij}=0$. Note that
$${\fontsize{10}{10}\selectfont\xy
(0,0)*{}; (12,12)*{} **\dir{-};
(3,0)*{}; (15,12)*{} **\dir{-};
(6,0)*{}; (18,12)*{} **\dir{-};
(9,0)*{}; (0,12)*{} **\dir{-};
(12,0)*{}; (3,12)*{} **\dir{-};
(15,0)*{}; (6,12)*{} **\dir{-};
(18,0)*{}; (9,12)*{} **\dir{-};
(0,12)*{}; (9,12)*{} **\dir{-};
(0,17)*{}; (9,17)*{} **\dir{-};
(0,12)*{}; (0,17)*{} **\dir{-};
(9,12)*{}; (9,17)*{} **\dir{-};
(12,12)*{}; (18,12)*{} **\dir{-};
(12,17)*{}; (18,17)*{} **\dir{-};
(12,12)*{}; (12,17)*{} **\dir{-};
(18,12)*{}; (18,17)*{} **\dir{-};
(4.5, 14)*{e_{i,n}};(0,-2)*{j}; (3,-2)*{j};(6,-2)*{j};
(15, 14)*{e_{j,m}};(9,-2)*{i}; (12,-2)*{i};(15,-2)*{i};(18,-2)*{i};(25,7)*{=};
\endxy}\ \ \ {\fontsize{10}{10}\selectfont\xy
(0,2)*{}; (12,14)*{} **\dir{-};
(3,2)*{}; (15,14)*{} **\dir{-};
(6,2)*{}; (18,14)*{} **\dir{-};
(9,2)*{}; (0,14)*{} **\dir{-};
(12,2)*{}; (3,14)*{} **\dir{-};
(15,2)*{}; (6,14)*{} **\dir{-};
(18,2)*{}; (9,14)*{} **\dir{-};
(0,-3)*{}; (6,-3)*{} **\dir{-};
(0,2)*{}; (6,2)*{} **\dir{-};
(0,-3)*{}; (0,2)*{} **\dir{-};
(6,-3)*{}; (6,2)*{} **\dir{-};
(9,-3)*{}; (18,-3)*{} **\dir{-};
(9,2)*{}; (18,2)*{} **\dir{-};
(9,-3)*{}; (9,2)*{} **\dir{-};
(18,-3)*{}; (18,2)*{} **\dir{-};
(3, -1)*{e_{j,m}};(0,16)*{i}; (3,16)*{i};(6,16)*{i};
(13.5, -1)*{e_{i,n}};(9,16)*{i}; (12,16)*{j};(15,16)*{j};(18,16)*{j};
\endxy}.$$
The right multiplication by
$${\fontsize{10}{10}\selectfont\xy
(0,0)*{}; (12,12)*{} **\dir{-};
(3,0)*{}; (15,12)*{} **\dir{-};
(6,0)*{}; (18,12)*{} **\dir{-};
(9,0)*{}; (0,12)*{} **\dir{-};
(12,0)*{}; (3,12)*{} **\dir{-};
(15,0)*{}; (6,12)*{} **\dir{-};
(18,0)*{}; (9,12)*{} **\dir{-};
(0,-2)*{j}; (3,-2)*{j};(6,-2)*{j};
(9,-2)*{i}; (12,-2)*{i};(15,-2)*{i};(18,-2)*{i};
\endxy}$$
is a map from $R(i^n j^m)\ac e_{i,n}\ot e_{j,m}$ to $R(i^n j^m)\ac e_{j,m}\ot e_{i,n}$, which has the obvious inverse by flipping this diagram.
\end{proof}
\end{proposition}

Let $K_0(R)_{\Q(q)}=\Q(q)\ot_{\Z[q,q^\zz]}K_0(R)$. By (\ref{P}) and Proposition \ref{Serre}, we have a well-defined bialgebra homomorphism
$$\begin{aligned}
& \Ga_{\Q(q)}: U^-\rightarrow K_0(R)_{\Q(q)}\\
&\phantom{\Ga_{\Q(q)}:a} F_i\ma [P_i]\qquad \tx{for}\ i\in I^\pm\\
& \phantom{\Ga_{\Q(q)}:} \ F_{in}\ma [P_{i,n}]\qquad \tx{for}\ i\in I^0, n>0.
\end{aligned}$$
By Proposition \ref{Sy}, the bilinear form $\{ \ , \ \}$ on $U^-$ and the Khovanov-Lauda's form $( \ , \ )$ on $K_0(R)_{\Q(q)}$ coincide under the map $\Ga_{\Q(q)}$, that is
$$(\Ga_{\Q(q)}(x),\Ga_{\Q(q)}(y))=\{x,y\} \ \ \tx{for} \ x,y\in U^-.$$
Thus, $\Ga_{\Q(q)}$ is injective by the non-degeneracy of $\{ \ , \ \}$. We have
$$ \Ga_{\Q(q)}(\ovv{x})=\ovv{\Ga_{\Q(q)}(x)}\ \ \tx{for} \ x\in U^-.$$
Moreover, $\Ga_{\Q(q)}$ induces an injective $\A$-bialgebra homomorphism
$\Ga:{_{\A}}U^- \ra K_0(R)$.

\subsection{Surjectivity of $\Ga_{\Q(q)}$ and $\Ga$}\

\vs{1mm}

Let $\nu\in \N[I]$. Define $\underline{\nu}$ to be the set of sequences $\ii$ of type $\nu$ with `parameters' in $i\in I^0$. Such a sequence is of the form
 $$\ii=j_1 \ds j_{p_0}  (i_1,n_1) k_1 \ds k_{p_1}(i_2,n_2)\ \ds\ (i_t,n_t) h_1 \ds h_{p_t},$$
with $(i_1,n_1),\ds,(i_t,n_t)\in I^0\ti \Z_{>0}$ and such that the expended sequence
 $$j_1 \ds j_{p_0}\underbrace{i_1\ds i_1}_{n_1}k_1 \ds k_{p_1}\underbrace{i_2\ds i_2}_{n_2}\ \ds\ \underbrace{i_t\ds i_t}_{n_t} h_1 \ds  h_{p_t}$$
belongs to $\Seq(\nu)$. For each $\ii\in \underline{\nu}$, we assign the following idempotent of $R(\nu)$
$$1_{\ii}=1_{j_1 \ds j_{p_0}}\ot e_{i_1,n_1}\ot 1_{k_1\ds k_{p_1}}\ot e_{i_2,n_2}\ot\ \cs\ \ot  e_{i_t,n_t}\ot 1_{h_1\ds h_{p_t}}.$$

Define the character of $M\in R(\nu)\fMod$ as
$$\Ch M=\s_{\ii\in\underline{\nu}}\ (\Dim 1_{\ii}M)\ii\ \in \Z[q,q^\zz]\underline{\nu}.$$
Note that, each $\ii\in \underline{\nu}$ determines a monomial $\Theta_{\ii}$ in $U^-$ under the correspondence
$$i\ma F_i\ \ \tx{for}\ i\in I^\pm, \quad (i,n)\ma F_{in}\ \ \text{for}\ i\in I^0,n>0.$$

Let $U^-_{\nu}$ be the $\Q(q)$-subspace of $U^-$ spanned by $\Theta_{\ii}$ for all $\ii\in \underline{\nu}$. Combining with $\Ga_{\Q(q)}$, we obtain a $\Q(q)$-linear map
$$\Q(q)\underline{\nu} \longrightarrow U^-_{\nu}\ov{\Ga_{\Q(q)}}{\longrightarrow}K_0(R(\nu))_{\Q(q)},$$
which has the dual map
$$G_0(R(\nu))_{\Q(q)}\ov{\Ch}{\longrightarrow}\Q(q)\underline{\nu}.$$
We next show that the character map $\Ch$ is injective.

Let $i\in I^0$ and $\nu=ni$. Then
$$\underline{\nu}=\{(i,n_1)\ds(i,n_t)\mid n_1+\ds + n_t=n\}.$$
Since $\Ga_{\Q(q)}: U^-_{ni}\ra K_0(R(\nu))_{\Q(q)}$ is injective and they have the same dimension $|\Pp_n|$, we see that $\Ga_{\Q(q)}$ is an isomorphism in this case, and so $\Ch: G_0(R(ni))_{\Q(q)}\rightarrow\Q(q)\underline{\nu}$ is injective.

\begin{lemma}\label{LLL1}
Let $i\in I^0$. The characters of all non-isomorphic gr-irreducible $R(ni)$-modules are $\Q(q)$-linearly independent.
\end{lemma}

\begin{example} Recall that the irreducible $\C[S_n]$-modules are {\it Specht modules} parametrized by the partitions $n$ (see, for example, \cite{CST2010}).

In this example, we deal with the case when $n=3$. The irreducible $\C[S_3]$-modules and their characters are given by
$$\begin{aligned}
& S^{(3)}=P_{i,3}=\C[S_3]\ac e_{i,3},\ \ \Ch \ S^{(3)}=e_{i,3}+e_{i,1}\ot e_{i,2}+e_{i,2}\ot  e_{i,1}+e_{i,1}\ot  e_{i,1}\ot  e_{i,1}\\
& S^{(21)}=\C[S_3]\ac 1/3(1+s_1-s_2s_1-s_1s_2s_1), \ \ \Ch\  S^{(21)}=e_{i,1}\ot e_{i,2}+e_{i,2}\ot  e_{i,1}+2e_{i,1}\ot e_{i,1}\ot e_{i,1}\\
& S^{(111)}=\C[S_3]\ac 1/6(1-s_1-s_2+s_1s_2+s_2s_1-s_1s_2s_1),\ \ \Ch \ S^{(111)}=e_{i,1}\ot  e_{i,1}\ot e_{i,1}
\end{aligned}.$$
Here, $S^{(3)},S^{(21)},S^{(111)}$ are the Specht modules.
\end{example}

Let  $i\in I$ and $n\geq 0$. Define a functor
\begin{equation*}
\begin{aligned}
\De_{i^n}: R(\nu)\Mod & \ra R(\nu-ni)\ot R(ni)\Mod \\
M & \ma (1_{\nu-ni}\ot 1_{ni})M.
\end{aligned}
\end{equation*}
For each $M\in R(\nu)\fMod$, we define
$$\vep_iM=\max\{n\geq 0\mid\De_{i^n}M\neq 0\}.$$
The following lemma can be proved by the same manner in \cite[Section 3.2]{KL2009}.

\begin{lemma}\label{LLL2}
Let $i\in I$ and  $M\in R(\nu)\fMod$ be a gr-irreducible module with $\vep_iM=n$. Then $\De_{i^n}M$ is isomorphic to $K\ot V$ for some gr-irreducible $K\in R(\nu-ni)\fMod$ with $\vep_iK=0$ and some gr-irreducible $V\in R(ni)\fMod$. Moreover, we have
$$M\co \hd\ \Ind_{\nu-ni,ni} K\ot V.$$
\end{lemma}

Recall that a gr-irreducible $R(ni)$-module $V$ has one of the following forms:

\begin{itemize}
\item[(i)] if $i \in I^+$, then $V=V(i^n)$, the unique gr-irreducible module of the nil-Hecke algebra,
\item[(ii)] if $i\in I^-$, then $V=V(i^n)$ is the one-dimensional trivial module,
\item[(iii)] if $i\in I^0$, then $V$ has $|\Pp_n|$ many choices.
\end{itemize}

\begin{theorem}
The map $\Ch: G_0(R(\nu))_{\Q(q)}\ra\Q(q)\underline{\nu}$ is injective.
\begin{proof}
We show that the characters of elements in $\B_{\nu}$ are linearly independent over $\Q(q)$ by induction on $\tx{ht}(\nu)$. The case of $\tx{ht}(\nu)=0$ is trivial. Assume for $\tx{ht}(\nu)<n$, our assertion is true. Now, suppose $\tx{ht}(\nu)=n$ and we are given a non-trivial linear combination
\begin{equation}\label{Ch}
\s_{M\in\B_{\nu}}c_M\Ch M=0
\end{equation}
for some $c_M\in\Q(q)$. Choose $i\in I$. We prove by a downward induction on $k=n,\ds,1$ that $c_M=0$ for all $M$ with $\vep_iM=k$.

If $k=n$ and $M\in \B_{\nu}$ with $\vep_iM=n$, then $\nu=ni$ and $M$ is a gr-irreducible $R(ni)$-module. When $i\in I^+\sqcup I^-$, our assertion is trivial. When $i\in I^0$, it follows from Lemma \ref{LLL1}.

Assume for $1\leq k<n$, we have $c_L=0$ for all $L$ with $\vep_iL>k$. Taking out the terms with $i^k$-tail in the rest of (\ref{Ch}), we obtain
\begin{equation}\label{Ch1}
\s_{M\colon\vep_iM=k}c_M\Ch (\De_{i^k}M)=0.
\end{equation}

By Lemma \ref{LLL2}, we can assume $\De_{i^k}M\co K_M\ot V_M$ for gr-irreducible $K_M\in R(\nu- ki)\fMod$  with $\vep_iK_M=0$ and gr-irreducible $V\in R(ki)\fMod$. Then (\ref{Ch1}) becomes $$\s_{M\colon\vep_iM=k}c_M\Ch K_M\ot\Ch V_M=0.$$

Note that if $[M]\neq [M']$ in $\B_{\nu}$, then we have $[K_M]\neq [K_{M'}]$ or $[V_M]\neq [V_{M'}]$. By the induction hypothesis, $\Ch K$ $(K\in\B_{\nu-ki})$ are linearly independent, and by  Lemma \ref{LLL1}, if $i\in I^0$, $\Ch V$ $(V\in\B_{ ki})$ are linearly independent. It follows that $c_M=0$ for all $M$ with $\vep_iM=k$.
Since each gr-irreducible $R(\nu)$-modules $M$ has $\vep_iM>0$ for at least one $i\in I$, the theorem has been proved.
\end{proof}
\end{theorem}

\begin{remark}
We see from the proof that the map $\Ch: G_0(R(\nu))\ra\Z[q,q^\zz]\underline{\nu}$ is also injective. If we set $\tx{ch}={\Ch |}_{q=1}$, then by a similar argument, the ungraded  characters of elements in $\B_{\nu}$ are linearly independent over $\Z$.
\end{remark}

\begin{corollary}
$\Ga_{\Q(q)}: U^-\ra K_0(R)_{\Q(q)}$ is an isomorphism.
\end{corollary}

We next consider the surjectivity of $\Ga:{_{\A}}U^-\ra K_0(R)$.

For $\la\vdash n$, let $S^\la$ denote the Specht module corresponding to $\la$ (or the Young diagram of shape $\la$). Let $\prec$ be the dominance ordering on $\Pp_n$.
For $\mu\prec \la$, we denote by $S^{\la/\mu}$ the shew representation of $\C[S_n]$ corresponding to the skew diagram $\la/\mu$.

\begin{lemma}\cite[Proposition 3.5.5]{CST2010}\label{LLL3}
Let $\la\vdash n$ be a partition and $(b_1,\ds,b_l) \vDash n$ be a composition of $n$. Then
$$\Res^n_{b_1,\ds,b_l}S^\la=\op(S^{\la^{(1)}}\ot S^{\la^{(2)}/\la^{(1)}}\ot \cs\ot S^{\la/\la^{(l-1)}}),$$
where the sum runs over all sequences $\la^{(1)}\prec \la^{(2)}\prec \cs\prec \la^{(l)}=\la$ such that  $|\la^{(j)}/\la^{(j-1)}|=b_j$ for all $j=1,\ds,l$.
\end{lemma}

\begin{lemma}\cite[Proposition 3.5.12]{CST2010}\label{LLL4}
 Assume $|\la/\mu|=k$. The multiplicity of the trivial representation $S^{(k)}$ in $S^{\la/\mu}$ is $1$ if $\la/\mu$ is totally disconnected, and $0$ otherwise.
\end{lemma}

\begin{proposition}\label{PPP1}
Let $\la\vdash n$. Let $\mathbf c \vDash n$ be a composition of $n$, which determines a partition $\la_{\mathbf c}\vdash n$. Then
$$\tx{dim}_\C(e_{i,\mathbf c}\ac S^\la)=
\begin{cases}1\quad \tx{if}\ \la_{\mathbf c}=\la,\\
0 \quad \tx{if}\ \la_{\mathbf c}>\la,\end{cases}$$
where $>$ is the lexicographical order of partitions.
\begin{proof}
The proposition follows from Lemma \ref{LLL3} and Lemma \ref{LLL4}, using the fact that the Kostka number
$$K_{\la, \mathbf c}=K_{\la, \la_{\mathbf c}}=
\begin{cases}1\quad \tx{if}\ \la_{\mathbf c}=\la,\\
0 \quad \tx{if}\ \la_{\mathbf c}>\la.\end{cases}$$
\end{proof}
\end{proposition}

Note that Lemma \ref{LLL1} can be derived directly from the proposition above.

For $i\in I^0,\la=(c_1,\ds,c_r)\vdash n$, we set
$$P_{i,\la}=P_{i,c_1}\cs P_{i,c_r}=R(ni)\ac e_{i,\la}.$$
Then $(P_{i,{\la}}, S^{\mu})=\tx{dim}_\C(e_{i,\la}\ac S^\mu)$ and according to Proposition \ref{PPP1}, the matrix $\left\{(P_{i,{\la}}, S^{\mu})\right\}_{\la,\mu \vdash n}$
unitriangular. It follows that each $[P]\in K_0(R(ni))$ can be written as a $\Z[q,q^\zz]$-linear combination of $[P_{i,\la}]$ for $\la \vdash n$.

More generally, we write the set $I$ as
$$I=\{i_1,i_2,\ds,i_k,\ds\}.$$
For a gr-irreducible $R(\nu)$-module $M$, we let $c_1=\vep_{i_1}M$ and assume  $\De{i_1^{c_1}}M = M_1\ot V_M$ for gr-irreducible $M_1$ with $\vep_{i_1} M_1=0$ and gr-irreducible $V_M\in R(c_1i_1)\Mod$.

If $i_1\in I^0$, then $V_M=S^{\la^{(1)}}$ for some $\la^{(1)}\vdash c_1$. So we get a pair $(c_1,\la^{(1)})$, where we set $\la^{(1)}=0$ when $i_1\notin I^0$. Inductively,  $c_k=\vep_{i_k}M_{k-1}$ and $\De_{{i_k}^{c_k}}M_{k-1} = M_k\ot V_{M_{k-1}}$, and we obtain $(c_k,\la^{(k)})$. If we do not get $M_k=0$ after $I$ exhausted, we can continue
the above process from $i_1$. Therefore, each $b\in \B_{\nu}$ is assigned by a sequence
$$W_b=(c_1,\la^{(1)})(c_2,\la^{(2)})\ds (c_k,\la^{(k)})\ds,$$
and we see from Lemma \ref{LLL2} that $b$ is uniquely determined by $W_b$.

Set
$$P_{(c_k,\la^{(k)})}=
\begin{cases} P_{{i_k}^{(c_k)}}&\tx{if}\ i_k \in I^+,\\
P_{{i_k}^{c_k}}& \tx{if}\ i_k\in I^-,\\
P_{{i_k}, \la^{(k)}} &\tx{if}\ i_k\in I^0,\end{cases}$$
and
$$P_{W_b}=\cs P_{(c_k,\la^{(k)})} \cs P_{(c_2,\la^{(2)})} P_{(c_1,\la^{(1)})}.$$

Let $b,b'\in \B_{\nu}$ with $W_b=(c_1,\la^{(1)})(c_2,\la^{(2)})\ds$ and $W_{b'}=(d_1,\mu^{(1)})(d_2,\mu^{(2)})\ds$. We denote $W_b>W_{b'}$ if for some $t$, $$(c_1,\la^{(1)})=(d_1,\mu^{(1)}),\ds,(c_{t-1},\la^{(t-1)})=(d_{t-1},\mu^{(t-1)})$$
but $c_t>d_t\ \text{or}\ c_t=d_t, \la^{(t)}>\mu^{(t)}$.

\begin{proposition}
$\HOM(P_{{W_b}}, S_{b'})=0$ if $b>b'$ and $\HOM(P_{{W_b}}, S_b)\co\C$.
\begin{proof}
For $i\in I^+$, we have $\HOM(P_{i^{(n)}},V(i^n))\co\C$ since $P_{i^{(n)}}$ is the graded projective cover of $V(i^n)$. For $i\in I^-$, $R(ni)=P_{i^n}$ is the graded projective cover of $V(i^n)$. The results follows immediately from the Frobenius reciprocity and Proposition \ref{PPP1}, which deals with the case when $i\in I^0$.
\end{proof}
\end{proposition}

By proposition above, each $[P]\in K_0(R(\nu))$ can be written as a $\Z[q,q^\zz]$-linear combination of $[P_{{W_b}}]$ for $b\in\B_{\nu}$. Therefore, $\Ga$ is surjective.

\begin{theorem}
$\Ga: {_\A}U^-\ra K_0(R)$ is an isomorphism.
\end{theorem}

For $M\in R(\nu)\fMod$, let $M^\divideontimes=\HOM_\C(M,\C)^\psi$ be the dual module in $R(\nu)\fMod$ with the action given by
$$(zf)(m):=f(\psi(z)m)\ \ \tx{for}\ z\in R(\nu), f\in \HOM_\C(M,\C), m\in M.$$
As proved in \cite[Section 3.2]{KL2009}, for each gr-irreducible $R(\nu)$-module $L$, there is a unique $r\in\Z$ such that $(L\{r\})^\divideontimes \co L\{r\}$, and the graded projective cover of $L\{r\}$ is stable under the bar-involution $^-$.

\section{\textbf{Geometric realization of $R$}}

In this chapter, we prove the following theorem:

\begin{theorem}\label{canonical}
Under $\Ga$, the canonical bases in $\bigsqcup_{\nu\in\N[I]}\PP_{\nu}$ of $\KK={_\A} U^-$ are mapped to the self-dual indecomposable projective modules in $K_0(R)$.
\end{theorem}

Let's consider the Jordan quiver first.

\subsection{Jordan quiver case}\

\vs{1mm}

Assume $I=I^0=\{i\}$ and ${_\A}U^-$ be the associated quantum group. Following the notations in  Remark \ref{iso}, it is already known that  the canonical basis of ${_\A}U^-_{ni}$ is the $ \{IC(\OO_{\la})\}_{\la\vdash n}$.

In the following lemma, we denote by $K_0(S_n)$ the Grothendieck group of the finite dimensional $\C[S_n]$-modules, then $\op_n K_0(S_n)$ is a $\Z$-algebra with the multiplication induced by the induction of modules.

\begin{lemma}\label{spt}
If we have an isomorphism $\Phi:\A\ot_\Z(\op_n K_0(S_n))\ra {_\A}U^-$ of $\A$-algebras, which sends the trivial representation $S^{(n)}$ to $F_{in}=IC(\OO_{(1^n)})$ or $IC(\OO_{(n)})$, then it maps the irreducible $\C[S_n]$-modules to the canonical bases of ${_\A}U^-$. More precisely:
\begin{itemize}
\item[(1)] If $\Phi(S^{(n)})=IC(\OO_{(1^n)})$, then $\Phi(S^\la)=IC(\OO_{\wi{\la}})$ for all $\la\vdash n$, \\
    where $\wi{\la}$ is the transpose of $\la$.
\item[(2)] If $\Phi(S^{(n)})=IC(\OO_{(n)})$, then $\Phi(S^\la)=IC(\OO_{\la})$ for all $\la\vdash n$.
\end{itemize}
\begin{proof}
Let $\La$ be the ring of symmetric functions. By \cite{Lus81} (see also \cite[Example 3.10]{S09}), there is an $\A$-algebra isomorphism
$$\A\ot_\Z\La\xs {_\A} U^-,$$
which sends the Schur functions $s_{\la}$ to $IC(\OO_{\la})$ for all $\la$. By the classical representation theory of the symmetric group, there is an $\A$-algebra isomorphism
$$\A\ot_\Z(\op_n K_0(S_n)) \xs \A\ot_\Z\La,$$
which sends the Specht module $S^\la$ to $s_{\la}$. Furthermore, we have an $\A$-algebra involution $\om$ of $\A\ot_{\Z}(\op_n K_0(S_n))$ sending $S^\la$ to $S^{\wi{\la}}$.
\end{proof}
\end{lemma}

Since $I=I^0=\{i\}$, we have $K_0(R)=\op_nK_0(R(ni))$. There is an obvious $\A$-algebra isomorphism
$$\Omega:\A\ot_\Z(\op_n K_0(S_n)) \xs K_0(R),$$
which sends the irreducible modules to their gr-projective covers. We thus obtain an $\A$-algebra isomorphism $\Omega'=\Ga^\zz\Omega:\A\ot_\Z(\op_n K_0(S_n))\ra {_\A}U^-$ with $\Omega'(S^{(n)})=F_{in}$. By Lemma \ref{spt}, we conclude that  $\Omega'(S^\la)=IC(\OO_{\wi{\la}})$ for all $\la\vdash n$. Therefore, $\Ga^\zz$ maps the self-dual indecomposable projective modules to the canonical bases.

\subsection{Springer functor}\

\vs{1mm}

For each $n$, we have the Springer functor (see e.g. \cite[Chapter 8]{Ach21})
$$\Hom(\Spr_n,-):P_{G_{ni}}(\EE_{ni}^{\tx{nil}})\ra \C[S_n]\tx{-mod}, $$
which is an equivalence of categories, mapping $\{IC(\OO_{\la})\}_{\la\vdash n}$ to irreducible modules. In particular,
\begin{equation}\label{sf}
\Hom(\Spr_n,IC(\OO_{(1^n)}))=S^{(n)},\quad \Hom(\Spr_n,\Spr_n)=\C[S_n],
\end{equation}

It has a quasi-inverse
$$\Spr_n\ot_{S_n}-: \C[S_n]\tx{-mod}\ra P_{G_{ni}}(\EE_{ni}^{\tx{nil}}).$$
According to \cite[Theorem 1.3]{C08}, this functor induces an algebra isomorphism
$$ \op_n\Spr_n\ot_{S_n}-: \A\ot_\Z(\op_n K_0(S_n)) \xs {_\A}U^-,$$
with the inverse $\op_{n}\Hom(\Spr_n,-)$. Thus, by Lemma \ref{spt}, we obtain $\Hom(\Spr_n, IC(\OO_{\la}))=S^{\wi{\la}}$ for all $\la\vdash n$.

\subsection{Proof of Theorem \ref{canonical}}\

\vs{1mm}

Assume we are given an arbitrary quiver $(I,H)$. Let $\nu\in\N[I]$. For $\ii,\jj\in \Seq(\nu)$, we define the Steinberg-type varieties
$$\ZZ_{\ii,\jj}=\wi{\FF}_{\ii}\ti_{\EE_\nu}\wi{\FF}_{\jj},\quad \ZZ_{\nu}=\bigsqcup_{\ii,\jj}\ZZ_{\ii,\jj},$$
where the fiber product is taken with respect to the projections $\pi_{\ii}$ and $\pi_{\jj}$.

Define the $G_{\nu}$-equivariant (Borel-Moore) homology spaces as
$$\R_{\ii,\jj}=H_{*}^{G_\nu}(\ZZ_{\ii,\jj}),\quad \R_\nu=\op_{\ii,\jj}\R_{\ii,\jj}.$$
Then $\R_\nu$ has a $\A$-algebra algebra structure defined by the convolution product (see \cite{VV2011} for details).
Let $\LL_\nu=\op_{\ii\in\Seq(\nu)}\LL_\ii$. It is known that
$$\R_\nu\co \Ext^*_{G_\nu}(\LL_{\nu},\LL_{\nu}).$$

For each $i\in I$, we set $P_i(u,v)= {(u-v)}^{h_i}$ and
$$\begin{aligned}
& P_i(u,v,w)=-\f{P_i(v,u)P_i(u,w)}{(u-v)(u-w)}-\f{P_i(u,w)P_i(v,w)}{(u-w)(v-w)}+\f{P_i(u,v)P_i(v,w)}{(u-v)(v-w)}\\
& P_i'(u,v,w)=\f{P_i(u,v)P_i(u,w)}{(u-v)(u-w)}+\f{P_i(u,w)P_i(w,v)}{(u-w)(v-w)}-\f{P_i(u,v)P_i(v,w)}{(u-v)(v-w)}.
\end{aligned}.$$
The following result was proved in \cite{KKP2013}, and a more general treatment can be found in \cite{Ju2013}.

\begin{proposition}\cite{KKP2013} \
$\R_\nu$ is generated by $1_{\ii}$, $x_{k,\ii}$ $(1\leq k\leq n)$, $\si_{k,\ii}$ $(1\leq k\leq n-1)$ for $\ii\in \Seq(\nu)$, such that $\R_{\ii,\jj}=1_\ii\R_\nu 1_\jj$, subject to the following defining relations: (we set $ x_k=\sum_{\ii}x_{k,\ii}$ and $\si_k=\sum_{\ii}\si_{k,\ii}$)
$$\begin{aligned}
& 1_{\ii} 1_{\jj} = \de_{\ii,\jj} 1_{\ii},\quad 1=\s_{\ii \in \Seq(\nu)} 1_{\ii}, \quad x_{k,\ii}=1_{\ii} x_k=x_k 1_\ii, \quad x_kx_l = x_lx_k,\\
& \si_{k,\ii}=1_{s_k{\ii}}\si_k=\si_k 1_\ii,\quad \si_k\si_l = \si_l\si_k\ \ \tx{if}\ |k-l|>1, \\
& \si_k^21_\ii=\begin{cases}-\pl_k \left(P_{i_k}(x_k,x_{k+1})\right)\ac \si_{k,\ii} & \tx{if}\  i_k = i_{k+1}, \\ Q_{i_k, i_{k+1}}(x_k,x_{k+1})1_\ii   & \tx{if}\ i_k \neq i_{k+1},
\end{cases} \\
&  (x_{s_k (l)}\si_k-\si_k x_{l})1_\ii=\begin{cases}
P_{i_k }(x_k, x_{k+1})1_\ii & \tx{if}\ l=k \ \tx{and}\ i_k = i_{k+1}, \\
-P_{i_k}(x_k, x_{k+1})1_\ii  & \tx{if}\ l=k+1 \ \tx{and}\ i_k = i_{k+1}, \\
  0 & \tx{otherwise}, \end{cases}\\
& (s_k \si_{k+1}\si_k - \si_{k+1}\si_k \si_{k+1})1_\ii \\
& \qquad \qquad =
\begin{cases}
P_{i_k}(x_k, x_{k+1},x_{k+2})\si_{k,\ii}+ P'_{i_k}(x_k, x_{k+1},x_{k+2})\si_{k+1,\ii}& \tx{if}\ i_k = i_{k+1} = i_{k+2},\\
-P_{i_k}(x_k,x_{k+2})\f{Q_{i_k,i_{k+1}}(x_k, x_{k+1})-Q_{i_k,i_{k+1}}(x_{k+2}, x_{k+1})}{x_k-x_{k+2}}1_\ii & \tx{if}\ i_k = i_{k+2}\neq i_{k+1}, \\
0 & \tx{otherwise}.
\end{cases}
\end{aligned}$$
\end{proposition}

By direct calculation, we obtain the following proposition.

\begin{proposition}{\label{inter}}
For each $\nu\in\N[I]$, there is an isomorphism $R(\nu)\xs \R_\nu$ of graded algebras, given by
$$1_\ii\ma 1_\ii, \quad x_{k,\ii}\ma x_{k,\ii}, \quad
\tau_{k,\ii}\ma \begin{cases} \si_{k,\ii}+(x_{k,\ii}-x_{k+1,\ii})^{-a_{i_ki_k}/2} & \tx{if}\ i_k=i_{k+1}\in I^{\leq 0},\\ -\si_{k,\ii} & \tx{if}\ i_k=i_{k+1}\in I^+,\\ \si_{k,\ii} & \tx{if}\ i_k\neq i_{k+1}.\end{cases}$$
\end{proposition}

By a similar argument in \cite[Section 4.6]{VV2011}, we have an additive functor
$Y_\nu:\QQ_\nu\ra \R_\nu\pMod$ given by
$$\LL\ma \Ext^*_{G_\nu}(\LL_\nu,\LL),$$
which is an equivalence of semisimple categories. It has a quasi-inverse given by
 $$\Pp\ma \LL_\nu\ot_{\R_\nu} \Pp.$$

Similar to \cite[Theorem 1.3]{C08}, the functors $Y_\nu$ $(\nu\in\N[I])$ induces a $\A$-algebra isomorphism
$$\Ga': ({_\A}U^- \co)\ \KK \xs \op_{\nu\in \N[I]} K_0(\R_\nu)\ (\co K_0(R)),$$
where the second isomorphism follows from Proposition \ref{inter}. Moreover, we have
\begin{equation}\label{ree}\Ga'(F_i^{(n)})= [P_{i^{(n)}}] \ \ \tx{for}\ i\in I^+,n>0,\quad \Ga'(F_{i})= [P_{i}]\ \ \tx{for}\ i\in I^{\leq 0}.\end{equation}

Let $i\in I^0, n>0$. Note that $\LL_{ni}$ coincides with the Springer sheaf $\Spr_n$. Thus, $$(\R_{ni})_0=\Ext^0_{G_{ni}}(\LL_{ni},\LL_{ni})\co \C[S_n],\quad \Ext^0_{G_{ni}}(\LL_{ni},-)=\Hom(\Spr_n,-).$$
By (\ref{sf}), the degree zero part of $Y_{ni}(IC(\OO_{(1^n)}))$ is the trivial $\C[S_n]$-module. Since $Y_{ni}(IC(\OO_{(1^n)}))$ is its graded projective cover, we have
\begin{equation}\label{imm}\Ga'(F_{in})=[P_{i,n}]\ \ \tx{for}\ i\in I^0,n>0.\end{equation}
Combining (\ref{ree}) and (\ref{imm}), we conclude that $\Ga'=\Ga$, which proves Theorem \ref{canonical}.

\section{\textbf{Categorification of irreducible highest weight modules}}

We show that the cyclotomic quotient $R^{\La}$ $(\La\in P^+)$ categorifies the irreducible highest weight module $V(\La)$ over $U$, which arises from $U^-$ through the Drinfeld double construction.

\subsection{The algebra $U$}\

\vs{1mm}

Given a Borcherds-Cartan matrix $A=(a_{ij})_{i,j\in I}$, we set
\begin{itemize}
\item[(1)] $P^\vee=(\op_{i\in I}\Z h_i)\oplus (\op_{i\in I}\Z d_i)$, called the {\it dual weight lattice},
\item[(2)] $\h=\Q\ot_\Z P^\vee$, called the {\it Cartan subalgebra},
\item[(3)] $P=\{\la\in \h^*\mid \la(P^\vee)\se\Z\}$, called the {\it weight lattice},
\item[(4)] $\Pi^\vee=\{h_i\in P^\vee\mid i\in I\}$, called the set of {\it simple coroots},
\item[(5)] $\Pi=\{\al_i\in P \mid i\in I\}$, called the set of {\it simple roots}, which is linearly independent over $\Q$ and satisfies
$$\al_j(h_i)=a_{ij},\ \al_j(d_i)=\de_{ij}\ \ \tx{for all}\ i,j \in I.$$
\end{itemize}

Denote by $P^+=\{\La\in P\mid \La(h_i)\geq 0 \ \tx{for all} \ i \in I\}$ the set of {\it dominant integral weights}. The free abelian group $Q=\op_{i \in I} {\Z \al_i}$ is called the {\it root lattice}. We identify $\Pi$ with $I$, and identify the positive root lattice $Q_+=\op_{i \in I}{\N\al_i}$ with $\N[I]$.

\begin{definition}
Define $\widehat{U}$ to be the $\Q(q)$-algebra generated by $q^h$ $(h\in P^{\vee})$ and $E_{in},
F_{in}$ $((i,n) \in \II)$, satisfying
\begin{align}
& q^0=1,\quad q^hq^{h'}=q^{h+h'} \quad \tx{for} \ h,h'\in P^{\vee}, \\
& q^hE_{jn}q^{-h}=q^{n\al_j(h)}E_{jn},\quad q^hF_{jn}q^{-h}=q^{-n\al_j(h)}F_{jn}\quad  \tx{for} \ h \in P^{\vee}, (j,n)\in\II, \\
& \s_{r+s=1-na_{ij}}(-1)^rE_i^{(r)}E_{jn}E_i^{(s)}=0 \quad \tx{for} \ i\in I^+, (j,n)\in \II \ \tx {and} \ i \neq (j,n), \\
& \s_{r+s=1-na_{ij}}(-1)^rF_i^{(r)}F_{jn}F_i^{(s)}=0 \quad  \tx{for} \ i\in I^+, (j,n)\in \II \ \tx {and} \ i \neq (j,n), \\
& E_{in}E_{jm}-E_{jm}E_{in}=F_{in}F_{jm}-F_{jm}F_{in}=0 \quad  \tx{for} \ a_{ij}=0.
\end{align}
We assign $|q^h|=0$, $|E_{in}|= n\al_i$ and $|F_{in}|= -n\al_i$ so that $\widehat{U}$ is $Q$-graded.
\end{definition}

The algebra $\widehat{U}$ is endowed with a co-multiplication
$\De:\widehat{U}\ra\widehat{U}\ot\widehat{U}$ given by
$$\De(E_{in}) = \s_{r+t=n} E_{ir}K_i^{-t}\ot E_{it}, \quad \De(F_{in}) =\s_{r+t=n} F_{ir}\ot K_i^{r}F_{it},\ \ \tx{for}\ (i,n)\in\II,$$
and $\De(q^h) = q^h\ot q^h$. Here $K_i=q^{h_i}$ for all $i\in I$.

Let $\widehat{U}^-$ (resp. $\widehat{U}^{\leq 0}$) be the subalgebra of $\widehat{U}$ generated by $F_{in}$ $((i,n)\in\II)$ (resp. $F_{in}$ $((i,n)\in\II)$ and $q^h$ $(h\in P^{\vee})$). We  identify $\widehat{U}^-$ with the algebra $U^-$ defined in Definition \ref{U-}, and extend the bilinear form $\{ \ , \ \}$ to $\widehat{U}^{\leq 0}$ by setting
$$\{q^h,1\}=1,\quad \{q^h,F_{in}\}=0,\quad \{q^h,K_j\}=q^{-\al_j(h)}.$$

\begin{definition} By the Drinfeld double process, we define the algebra $U$ as the quotient of $\widehat{U}$ by the relations
\begin{equation}\label{drinfeld}
\s\{a_{(1)},b_{(2)}\}\om(b_{(1)})a_{(2)}=\s\{a_{(2)},b_{(1)}\}a_{(1)}\om(b_{(2)})\ \ \tx{for all}\ a,b \in \widehat{U}^{\leq 0}.
\end{equation}
Here we use the Sweedler's notation
$\De(x)=\sum x_{(1)}\ot x_{(2)}$, and $\om$ is the $\Q(q)$-algebra involution of $\widehat{U}$ given by $\om(E_{in})=F_{in}$, $\om(F_{in})=E_{in}$ and $\om(q^h)=q^{-h}$.
\end{definition}

According to \cite[Appendix]{FKKT2021}, one can obtain from equation (\ref{drinfeld}) that
\begin{equation} \label{comm}
[E_{jl},F_{it}]=
\begin{cases}
0  \ \ & \tx{if}\ i\neq j, \\
\sum_{p=1}^{\min\{l,t\}}\ga_pF_{i,t-p}E_{i,l-p} \ \ & \tx{if}\ i=j,
\end{cases}
\end{equation}
where the elements $\{\ga_p\}_{p\geq 1}$ are defined inductively:
$$\ga_p=\ep_p(K_i^{-p}-K_i^p)-\s_{r=1}^{p-1}\ep_rK_i^r\ga_{p-r}.$$
We will show in the Lemma \ref{coeff} that
$$\ga_p=q^{-p}\p_{k=1}^p\f{q^{k-1}K_i-q^{1-k}K_i^\zz}{q^k-q^{-k}}.$$

Let $\La\in P^+$. The irreducible highest weight $U$-module $V(\La)$ is given by
$$\begin{aligned}
& V(\La)\co U\bigg/(\s_{(i,n)\in \II}UE_{in}+\s_{h\in P^\vee}U(q^h-q^{\La(h)})+\s_{i\in I^+}U F_i^{\La(h_i)+1}+\s_{i\in I^{\leq 0} \ \tx{with}\ \La(h_i)=0}UF_{in}) \\
& \phantom{V(\La}\ \co U^-\bigg/(\s_{i\in I^+}U^- F_i^{\La(h_i)+1}+\s_{i\in I^{\leq 0} \ \tx{with}\ \La(h_i)=0}U^-F_{in}),
\end{aligned}$$
where $U^-$ is the subalgebra of $U$ generated by $F_{in}$ $((i,n)\in \II)$, which can be identified with the Definition \ref{U-}.

\begin{definition}\label{int}
A $U$-module $M$ is said to be {\it integrable} if it satisfies
\begin{itemize}
\item [(1)] $M=\op_{\mu \in P}M_{\mu}$ with each $\tx{dim}M_{\mu}<\infty$,
where $$M_{\mu}=\{m\in M\mid q^{h}m=q^{\mu(h)}m\ \tx{for all}\ h\in P^{\vee}\}.$$
\item [(2)] $\wt(M):=\{\mu\in P\mid M_{\mu}\neq 0\}\se \bigcup_{i=1}^{s}(\la_i-Q_+)$ for finitely many $\la_1,\ds,\la_s\in P$.
\item [(3)] $E_i$ and $F_i$ are locally nilpotent on $M$ for each $i\in I^+$.
\item [(4)] If $i\in I^{\leq 0}$, then $\mu(h_i)\in\Z_{\geq 0}$ for all $\mu \in \wt(M)$.
\item [(5)] If $i\in I^{\leq 0}$ with $\mu(h_i)=0$, then $F_{in}M_{\mu}=0$.
\end{itemize}
\end{definition}

It is well known that an integrable highest weight module $M$ with highest weight $\La\in P^+$ is isomorphic to $V(\La)$ (see, for example, \cite[Corollary 5.7]{KK2020}).

Let $\ovv U$ be the algebra with the same generators and the defining relations as $U$ except the quantum Serre relations (4.3) and (4.4). Similar to \cite[Proposition B1]{KMPY1996}, we have the following lemma.

\begin{lemma}\label{Int}
Let $M$ be a $\ovv U$-module satisfying (1)-(5) in Definition \ref{int}. Then $M$ admits a natural structure of a $U$-module.
\begin{proof}
Let $S$ be the antipode of $\ovv U$. Then $\End_\C(M)$ carries a natural $\ovv U$-module structure given by
$$(x\ac f)(m)=\s x_{(1)}f(S(x_{(2)})m),\ \ \tx{where}\ \De(x)=\s x_{(1)}\ot x_{(2)}.$$
We also endow $\ovv U$ itself with a $\ovv U$-module structure via the adjoint action
$$x\ac u=\ad_x(u):=\s x_{(1)}uS(x_{(2)}).$$
With respect to these module structures, the homomorphism $X:\ovv{U}\ra\End(M)$ induced by the $\ovv U$-module structure on $M$ is $\ovv{U}$-lin.

Let $i\in I^+,(j,n)\in \II$ with $i\neq(j,n)$. Then we have
$$F_iX(E_{jn})=X(\ad_{F_i}(E_{jn}))=0,\quad K_i(E_{jn})=X(\ad_{K_i}(E_{jn}))=q^{na_{ij}}X(E_{jn}).$$
Since $\End(M)$ is integral as $U_{(i)}$-module, its semisimplicity implies $E_i^{1-na_{ij}}X(E_{jn})=0$.

Now, for $i\in I^+$ and $r\in\N$, we recall that
$$\De(E_i^r)=\s_{k=0}^rq^{k(r-k)}{\begin{bmatrix} r \\ k \end{bmatrix}} K_i^{-k} E_i^{r-k}\ot E_i^k,\quad \De(F_i^r)=\s_{k=0}^rq^{k(r-k)}{\begin{bmatrix} r \\ k \end{bmatrix}} F_i^k\ot K_i^kF_i^{r-k}.$$
$$S(E_i^r)=(-1)^rq^{-r(r-1)}K_i^rE_i^r,\quad S(F_i^r)=(-1)^rq^{r(r-1)}F_i^r K_i^{-r},\quad S(q^h)=q^{-h}$$
Consequently, the adjoint action of $E_i^r$ on $u\in \ovv U$ is given by
$$\begin{aligned}
& \ad_{E_i^r}(u)=\s_{k=0}^r(-1)^kq^{k(r+1-2k)}{\begin{bmatrix} r \\ k \end{bmatrix}} K_i^{-k} E_i^{r-k}uK_i^kE_i^k\\
& {\phantom{\ad_{E_i^r}(u)}}=\s_{k=0}^r(-1)^kq^{k(1-r)}{\begin{bmatrix} r \\ k \end{bmatrix}} E_i^{r-k}K_i^{-k} uK_i^kE_i^k.
\end{aligned}$$
Therefore, the relation $E_i^{1-na_{ij}}X(E_{jn})=0$ implies the quantum Serre relation (4.3). Applying the involution $\om$ of $\ovv U$, we obtain the relation (4.4) as well.
\end{proof}
\end{lemma}

\subsection{Functors $\E_{in},\F_{in}$ and $\ovv\F_{in}$}\

\vs{1mm}

Let $R(\nu)$ $(\nu\in\N[I])$ be the KLR algebra defined in Section 2.2. For $n\in\N$, set $$R(n)=\op_{\nu\colon\tx{ht}(\nu)=n}R(\nu).$$
We will denote by $x_k$ the element
$\sum_{\ii\in\Seq(\nu)}x_{k,\ii}$ in $R(\nu)$ (or $\sum_\nu\sum_{\ii\in\Seq(\nu)}x_{k,\ii}$ in $R(n)$); an analogous convention is adopted for $\tau_k$. By abuse of notation, for $\om\in S_n$ we sometimes write $\om$ instead of $\tau_\om$, where $\tau_\om$ is determined by a fixed reduced expression of $\om$.

For any $(i,n)\in\II$ and $\nu\in\N[I]$, define the functors
$$\begin{aligned}
& \E_{in}: R(\nu+ni)\Mod\ra R(\nu)\Mod,\  N\ma 1_\nu\ot e_{i,n}\hs{1mm} N=1_\nu\ot e_{i,n}\hs{1mm} R(\nu+ni)\ot_{R(\nu+ni)}N,\\
& \F_{in}: R(\nu)\Mod\ra R(\nu+ni)\Mod,\  M\ma (R(\nu+ni)\hs{1mm} 1_\nu\ot e_{i,n})\ot_{R(\nu)}M,\\
& \ovv\F_{in}: R(\nu)\Mod\ra R(\nu+ni)\Mod, \  M\ma (R(\nu+ni)\hs{1mm} e_{i,n}\ot 1_\nu)\ot_{R(\nu)}M.
\end{aligned}$$
Here $e_i=e_{i,1}=1_i$ when $i\in I^\pm$. 

We also identify the functors $\E_{in},\F_{in}$ and $\ovv\F_{in}$ with their defining kernels:
$$\E_{in}=1_\nu\ot e_{i,n}\hs{1mm} R(\nu+ni),\quad \F_{in}=R(\nu+ni)\hs{1mm} 1_\nu\ot e_{i,n},\quad
\ovv\F_{in}=R(\nu+ni)\hs{1mm} e_{i,n}\ot 1_\nu.$$
In particular, $\E_{in}$ is an $R(\nu)\tx{-}R(\nu+ni)$-bimodule, while $\F_{in}$ and $\ovv\F_{in}$ are $R(\nu+ni)\tx{-}R(\nu)$-bimodules.

For any $m,n\in\N$, let $D_{m,n}$ (resp. $D^\zz_{m,n}$) denote the set of minimal-length left (resp. right) $S_m\ti S_n$-coset representatives in $S_{m+n}$. Note that if $i\in I^0$, then for any $k\leq n$ we have
$$e_{i,n}=\f{k!(n-k)!}{n!}\big(\s_{v\in D_{k,n-k}}\hs{-1mm} v\ \big)\x e_{i,k}\ot e_{i,n-k}=\f{k!(n-k)!}{n!}e_{i,k}\ot e_{i,n-k}\x \big(\s_{u\in D^\zz_{k,n-k}}\hs{-1mm} u\ \big).$$

\begin{proposition}\label{CC1} We have the following natural isomorphisms:
\begin{itemize}
\item [(1)] If $i\in I^0$, then for any $l,t\geq 1$, $$\E_{il}\F_{it}\simeq\op_{p=0}^{\min\{l,t\}}\F_{i,t-p}\E_{i,l-p}\ot Z_p: R(\nu)\Mod\ra R(\nu+(t-l)i)\Mod,$$
  where $Z_p$ is the algebra of symmetric polynomials in $x_1,\ds,x_p$, each of degree $2$.
\item [(2)] If $i\in I^\pm$, then $$\E_i\F_i\simeq q^{-a_{ii}}\F_i\E_i\oplus(\tx{Id}\ot\C[t_i]): R(\nu)\Mod\ra R(\nu)\Mod,$$
    where $t_i$ is an indeterminate of degree $2$.
\item [(3)] If $(j,l),(i,t)\in\II$ with $i\neq j$, then
$$\E_{jl}\F_{it}\simeq q^{-lt a_{ij}}\F_{it}\E_{jl}: R(\nu)\Mod\ra R(\nu+ti-lj)\Mod.$$
\end{itemize}
\begin{proof} 
First assume that $i\in I^0$. Let $\nu\in\N[I]$ with $\tx{ht}(\nu)=n$. Without loss of generality, we may assume that $t\geq l$, and set $\nu'=\nu+(t-l)i$.
By the PBW-type basis \eqref{PBW} of the KLR algebra and the double coset decomposition of permutations, we obtain the decomposition (see the diagram below):
$$1_{\nu',li}R(\nu+ti)1_{\nu,ti}=\op_{k=0}^l\ \op_{\substack{v'\in D_{n-k,t-l+k}\\ v\in D_{k,l-k}}} 1_{\nu'}\ot 1_{li}\x v'\ot v \x v_k \x R(\nu)\ot R(ti).$$

$${\fontsize{10}{10}\selectfont
\xy
(-6,0)*{}; (-6,12)*{} **\dir{-};(-6,0)*{}; (6,0)*{} **\dir{-};(-6,-9)*{}; (6,-9)*{} **\dir{-};(-6,0)*{}; (-6,-9)*{} **\dir{-};(6,0)*{}; (6,-9)*{} **\dir{-}; (0,-5)*{R(\nu)};
(10,0)*{}; (28,0)*{} **\dir{-};(10,-9)*{}; (28,-9)*{} **\dir{-};(10,0)*{}; (10,-9)*{} **\dir{-};(28,0)*{}; (28,-9)*{} **\dir{-}; (19,-5)*{R(ti)};
(-6,12)*{}; (0,20)*{} **\dir{-};(-3,12)*{}; (9,20)*{} **\dir{-};
(0,12)*{}; (-6,20)*{} **\dir{-};(3,12)*{}; (-3,20)*{} **\dir{-};(6,12)*{}; (3,20)*{} **\dir{-};(9,12)*{}; (6,20)*{} **\dir{-};
(13,12)*{}; (16,20)*{} **\dir{-};(16,12)*{}; (22,20)*{} **\dir{-};(19,12)*{}; (28,20)*{} **\dir{-};(22,12)*{}; (13,20)*{} **\dir{-};(25,12)*{}; (19,20)*{} **\dir{-};(28,12)*{}; (25,20)*{} **\dir{-};
(-3,0)*{}; (-3,12)*{} **\dir{-};
(0,0)*{}; (13,12)*{} **\dir{-};
(3,0)*{}; (16,12)*{} **\dir{-};
(6,0)*{}; (19,12)*{} **\dir{-};
(10,0)*{}; (0,12)*{} **\dir{-};
(13,0)*{}; (3,12)*{} **\dir{-};
(16,0)*{}; (6,12)*{} **\dir{-};
(19,0)*{}; (9,12)*{} **\dir{-};
(22,0)*{}; (22,12)*{} **\dir{-};
(25,0)*{}; (25,12)*{} **\dir{-};
(28,0)*{}; (28,12)*{} **\dir{-};(31,15.5)*{v};(-9,15.5)*{v'};(47,6)*{v_k\in D_{n,t}\cap D^\zz_{n+t-l,l}};
(16,10)*{\underbrace{ \quad \  }_{\ \ k}};
\endxy}$$
Thus,
$$\E_{il}\F_{it}=1_{\nu'}\ot{e_{i,l}}\hs{1mm}R(\nu+ti)\hs{1mm}1_\nu\ot e_{i,t}=\op_{k=0}^l\ \op_{v'\in D_{n-k,t-l+k}} v'\ot e_{i,l} \x v_k \x R(\nu)\ot R(ti)e_{i,t}.$$

On the other hand,
$$\begin{aligned}
& \F_{i,t-l+k}\E_{ik}=(R(\nu')\hs{1mm} 1_{\nu-ki}\ot e_{i,t-l+k})\ot_{\nu-ki}(1_{\nu-ki}\ot e_{i,k}\hs{1mm} R(\nu))\\
& \phantom{\F_{i,t-l+k}\E_{ik}}=\op_{v'\in D_{n-k,t-l+k}} \Big(v'\x 1_{\nu-ki}\ot R((t-l+k)i) e_{i,t-l+k}\Big)\ot_{\nu-ki}\Big(1_{\nu-ki}\ot e_{i,k}\hs{1mm}R(\nu)\Big).
\end{aligned}$$
Let $\wi v=\sum_{v\in D_{k,l-k}}v$ and $\wi u=\sum_{u\in D^\zz_{t-l+k,l-k}}u$. Note that
$$\begin{aligned}
& (1_{\nu'}\ot\wi v)\x \Big(v'\x 1_{\nu-ki}\ot R((t-l+k)i) e_{i,t-l+k}\Big)\x v_k\ot \hs{.3mm} e_{i,l-k}R((l-k)i)e_{i,l-k}\x\Big(1_{\nu-ki}\ot e_{i,k}\hs{1mm} R(\nu)\Big)\x(1_\nu\ot\wi u)\\
& =(1_{\nu'}\ot \wi v)\x (v'\ot e_{i,k}\ot e_{i,l-k})\x v_k\x R(\nu)\ot R((t-l+k)i) e_{i,t-l+k}\ot R((l-k)i) e_{i,l-k}\x (1_\nu\ot \wi u)\\
& =v'\ot e_{i,l}\x v_k\x R(\nu)\ot R(ti)e_{i,t}.
\end{aligned}$$

Hence for each $f\in e_{i,l-k}R((l-k)i)e_{i,l-k}$, the assignment
$$x\ot y\ \ma\ (1_{\nu'}\ot\wi v) \x x \x (v_k\ot f)\x y \x (1_\nu\ot \wi u)$$
defines an injective homomorphism from $\F_{i,t-l+k}\E_{ik}$ to $\E_{il}\F_{it}$. Since $e_{i,l-k}R((l-k)i)e_{i,l-k}\co Z_{l-k}$, we obtain
$$\E_{il}\F_{it}\simeq\op_{k=0}^l\F_{i,t-l+k}\E_{ik}\ot Z_{l-k}\ \ \ov{p=l-k}{\rule{2em}{0.4pt}\kern-2em\raise0.6ex\hbox{\rule{2em}{0.4pt}}}\ \ \op_{k=0}^l\F_{i,t-p}\E_{i,l-p}\ot Z_p.$$

The case where $i\in I^-$ follows by the same proof. When $i\in I^+$, the statement was established in \cite[Theorem 3.6]{KK2012}.

Now let $(j,l),(i,t)\in\II$ with $i\neq j$. Assume $t\geq l$ and set $\nu'=\nu+ti-lj$. Since $i\neq j$, we have the decomposition
$$1_{\nu',lj}R(\nu+ti)1_{\nu,ti}=\op_{v'\in D_{n-l,t}}  v'\ot 1_{lj} \x v_l\x R(\nu)\ot R(ti),$$
and therefore
$$\E_{jl}\F_{it}=\op_{v'\in D_{n-l,t}} v'\ot e_{j,l}\x  v_l\x  R(\nu)\ot R(ti)e_{i,t}.$$

$${\fontsize{10}{10}\selectfont
\xy
(-6,0)*{}; (-6,12)*{} **\dir{-};(-6,0)*{}; (6,0)*{} **\dir{-};(-6,-9)*{}; (6,-9)*{} **\dir{-};(-6,0)*{}; (-6,-9)*{} **\dir{-};(6,0)*{}; (6,-9)*{} **\dir{-}; (0,-5)*{R(\nu)};
(10,0)*{}; (19,0)*{} **\dir{-};(10,-9)*{}; (19,-9)*{} **\dir{-};(10,0)*{}; (10,-9)*{} **\dir{-};(19,0)*{}; (19,-9)*{} **\dir{-}; (14.5,-5)*{R(ti)};
(-6,12)*{}; (0,20)*{} **\dir{-};(-3,12)*{}; (9,20)*{} **\dir{-};
(0,12)*{}; (-6,20)*{} **\dir{-};(3,12)*{}; (-3,20)*{} **\dir{-};(6,12)*{}; (3,20)*{} **\dir{-};(9,12)*{}; (6,20)*{} **\dir{-};
(-3,0)*{}; (-3,12)*{} **\dir{-};
(0,0)*{}; (13,12)*{} **\dir{-};
(3,0)*{}; (16,12)*{} **\dir{-};
(6,0)*{}; (19,12)*{} **\dir{-};
(10,0)*{}; (0,12)*{} **\dir{-};
(13,0)*{}; (3,12)*{} **\dir{-};
(16,0)*{}; (6,12)*{} **\dir{-};
(19,0)*{}; (9,12)*{} **\dir{-};
(-9,15.5)*{v'};
(16,15)*{\overbrace{ \quad \  }^{l}};
\endxy}$$

On the other hand,
$$\begin{aligned}
& \F_{it}\E_{jl}=(R(\nu')\hs{1mm} 1_{\nu-l j}\ot e_{i,t})\ot_{\nu-l j}(1_{\nu-l j}\ot e_{j,l}\hs{1mm} R(\nu))\\
&\phantom{\F_{it}\E_{jl}}=\op_{v'\in D_{n-l,t}} \Big(v'\x 1_{\nu-lj}\ot R(ti) e_{i,t}\Big)\ot_{\nu-lj}\Big(1_{\nu-lj}\ot e_{j,l}\hs{1mm} R(\nu)\Big).
\end{aligned}$$
Since
$$\Big(v'\x 1_{\nu-lj}\ot R(ti) e_{i,t}\Big)\x v_l\x\Big(1_{\nu-lj}\ot e_{j,l}\hs{1mm}R(\nu)\Big)=v'\ot e_{j,l}\x  v_l\x R(\nu)\ot R(ti)e_{i,t},$$
and $|v_l|=-lta_{ij}$, we see that $\E_{jl}\F_{it}\co q^{-lt a_{ij}}\F_{it}\E_{jl}$ as  $R(\nu')\tx{-}R(\nu)$-bimodules.
\end{proof}
\end{proposition}

\begin{proposition}\  
\begin{itemize}
\item [(1)] If $i\in I^0$, then for any $l,t\geq 1$, we have the natural isomorphism $$\E_{il}\ovv\F_{it}\simeq\op_{p=0}^{\min\{l,t\}}q^{-p(i\ac \nu)}\ovv\F_{i,t-p}\E_{i,l-p}\ot Z_p: R(\nu)\Mod\ra R(\nu+(t-l)i)\Mod.$$
\item [(2)] If $i\in I^-$, then $\E_i\ovv\F_i\simeq q^{-i\ac\nu}\ovv\F_i\E_i\oplus(\tx{Id}\ot\C[t_i])$.
\item [(3)] If $i\in I^+$, then there is a short exact sequence of $R(\nu)\tx{-}R(\nu)$-bimodules:
$$0\ra \ovv\F_i\E_i \ra \E_i\ovv\F_i\ra q^{-i\ac\nu}R(\nu)\ot\C[t_i]\ra 0.$$
\item [(4)] If $(j,l),(i,t)\in\II$ with $i\neq j$, then $\E_{jl}\ovv\F_{it}\simeq \ovv\F_{it}\E_{jl}$.
\end{itemize}
\begin{proof}
For parts (1) and (2), the proof is similar to that of Proposition \ref{CC1}, by considering the diagram below; we omit the details.

$${\fontsize{10}{10}\selectfont
\xy
(-6,0)*{}; (-6,12)*{} **\dir{-};(-6,0)*{}; (6,0)*{} **\dir{-};(-6,-9)*{}; (6,-9)*{} **\dir{-};(-6,0)*{}; (-6,-9)*{} **\dir{-};(6,0)*{}; (6,-9)*{} **\dir{-}; (0,-5)*{R(ti)};
(10,0)*{}; (28,0)*{} **\dir{-};(10,-9)*{}; (28,-9)*{} **\dir{-};(10,0)*{}; (10,-9)*{} **\dir{-};(28,0)*{}; (28,-9)*{} **\dir{-}; (19,-5)*{R(\nu)};
(-6,12)*{}; (0,20)*{} **\dir{-};(-3,12)*{}; (9,20)*{} **\dir{-};
(0,12)*{}; (-6,20)*{} **\dir{-};(3,12)*{}; (-3,20)*{} **\dir{-};(6,12)*{}; (3,20)*{} **\dir{-};(9,12)*{}; (6,20)*{} **\dir{-};
(13,12)*{}; (16,20)*{} **\dir{-};(16,12)*{}; (22,20)*{} **\dir{-};(19,12)*{}; (28,20)*{} **\dir{-};(22,12)*{}; (13,20)*{} **\dir{-};(25,12)*{}; (19,20)*{} **\dir{-};(28,12)*{}; (25,20)*{} **\dir{-};
(-3,0)*{}; (-3,12)*{} **\dir{-};
(0,0)*{}; (13,12)*{} **\dir{-};
(3,0)*{}; (16,12)*{} **\dir{-};
(6,0)*{}; (19,12)*{} **\dir{-};
(10,0)*{}; (0,12)*{} **\dir{-};
(13,0)*{}; (3,12)*{} **\dir{-};
(16,0)*{}; (6,12)*{} **\dir{-};
(19,0)*{}; (9,12)*{} **\dir{-};
(22,0)*{}; (22,12)*{} **\dir{-};
(25,0)*{}; (25,12)*{} **\dir{-};
(28,0)*{}; (28,12)*{} **\dir{-};(31,15.5)*{v};(-9,15.5)*{v'};(32,6)*{v_k};
(16,10)*{\underbrace{ \quad \  }_{\ \ k}};
\endxy}$$

\vs{2mm}

When $i\in I^+$, part (3) was established in \cite[Theorem 3.9]{KK2012}.

Let $(j,l),(i,t)\in\II$ with $i\neq j$. Assume $t\geq l$ and set $\nu'=\nu+ti-lj$. We have the decomposition
$$1_{\nu',lj}R(\nu+ti)1_{ti,\nu}=\op_{v\in D_{t,n-l}}  v\ot 1_{lj} \x R(ti)\ot R(\nu),$$
and therefore
$$\E_{jl}\ovv F_{it}=\op_{v\in D_{t,n-l}} v\ot e_{j,l} \x  R(ti)e_{i,t}\ot R(\nu).$$

$${\fontsize{10}{10}\selectfont\xy
(-6,0)*{}; (-6,12)*{} **\dir{-};(-6,0)*{}; (6,0)*{} **\dir{-};(-6,-9)*{}; (6,-9)*{} **\dir{-};(-6,0)*{}; (-6,-9)*{} **\dir{-};(6,0)*{}; (6,-9)*{} **\dir{-}; (0,-5)*{R(ti)};
(10,0)*{}; (28,0)*{} **\dir{-};(10,-9)*{}; (28,-9)*{} **\dir{-};(10,0)*{}; (10,-9)*{} **\dir{-};(28,0)*{}; (28,-9)*{} **\dir{-}; (19,-5)*{R(\nu)};
(-6,12)*{}; (-2,20)*{} **\dir{-}; (-3,12)*{}; (1,20)*{} **\dir{-};(0,12)*{}; (7,20)*{} **\dir{-};(3,12)*{}; (10,20)*{} **\dir{-};(6,12)*{}; (13,20)*{} **\dir{-}; (10,12)*{}; (-6,20)*{} **\dir{-};(13,12)*{}; (3,20)*{} **\dir{-}; (-6,12)*{}; (13,12)*{} **\dir{-};
(-3,0)*{}; (-3,12)*{} **\dir{-};
(0,0)*{}; (0,12)*{} **\dir{-};
(3,0)*{}; (3,12)*{} **\dir{-};
(6,0)*{}; (6,12)*{} **\dir{-};
(10,0)*{}; (10,12)*{} **\dir{-};
(13,0)*{}; (13,12)*{} **\dir{-};
(16,0)*{}; (16,12)*{} **\dir{-};
(19,0)*{}; (19,12)*{} **\dir{-};
(22,0)*{}; (22,12)*{} **\dir{-};
(25,0)*{}; (25,12)*{} **\dir{-};
(28,0)*{}; (28,12)*{} **\dir{-};
(22,15)*{\overbrace{ \qquad \quad \ }^{l}}; (-8,15.5)*{v}
\endxy}$$

\vs{2mm}

On the other hand,
$$\begin{aligned}
&\ovv\F_{it}\E_{jl}=(R(\nu')\hs{1mm}e_{i,t}\ot 1_{\nu-lj})\ot_{\nu-lj}(1_{\nu-lj}\ot e_{j,l}\hs{1mm}R(\nu)) \\
& \phantom{\ovv\F_{it}\E_{jl}}=\op_{v\in D_{t, n-l}} \Big(v\x R(ti)\hs{1mm} e_{i,t}\ot 1_{\nu-lj}\Big)\ot_{\nu-lj}\Big(1_{\nu-lj}\ot e_{j,l}\hs{1mm}R(\nu)\Big).
\end{aligned}$$
Thus $\E_{jl}\ovv\F_{it}\simeq \ovv\F_{it}\E_{jl}$.
\end{proof}
\end{proposition}

\subsection{Cyclotomic KLR algebras}\

\vs{1mm}

For $\La\in P^+$ and $i\in I$, set $\La_i=\La(h_i)\geq 0$. For each $\nu\in\N[I]$, the cyclotomic KLR algebra $R^{\La}(\nu)$ is defined by
$$R^\La(\nu)=\f{R(\nu)}{R(\nu)a^\La(x_1)R(\nu)},\ \ \tx{where}\ a^\La(x_1)=\s_{\ii\in\Seq(\nu)}x_{1,\ii}^{\La_{i_1}}.$$
That is, $R^{\La}(\nu)$ the quotient of $R(\nu)$ by the two sided ideal generated by $x_{1,\ii}^{\La_{i_1}}$ for all $\ii\in \Seq(\nu)$.

The following lemma implies the finiteness of cyclotomic KLR algebras. The proof is similar to \cite[Proposition 2.3]{LV2009}.

\begin{lemma}\
\begin{itemize}
\item [(1)] For $i\in I^+\cup I^0$, if we impose $x_1^b=0$ in $R(2i)$ for some $b\in \Z_{>0}$, then $x_2^b=0$.
\item [(2)] For $i\in I^-$, let $a=-a_{ii}$, then $x_1^b=0$ implies $x_2^{b+ab}=0$ in $R(2i)$.
\item [(3)] For $i\neq j$, if $a_{ij}=0$, then $x_{1,ji}^b=0$ implies $x_{2,ij}^b=0$ in $R(i+j)$;\\
if $a=-a_{ij}>0$, then $x_{1,ij}^b=x_{1,ji}^{b'}=0$  implies $x_{2,ij}^{b'+ab}=0$.
\end{itemize}
\begin{proof}
The case $i\in I^0$ is immediate. For $i\in I^+$, the statement is well known for the cyclotomic nil-Hecke algebra. We prove only (2), since the proof of (3) is similar.

Assume $i\in I^-$ and set $a=-a_{ii}$. Since $\tau_1^2=(-1)^{a/2}(x_1-x_2)^{a}$, we see that

$$0=\ {\fontsize{10}{10}\selectfont\xy
(0,5)*{}; (0,-5)*{} **\crv{(8,0)};(5,5)*{}; (5,-5)*{} **\crv{(-3,0)}?(.1)*{\bu};(0,-7)*{i}; (5,-7)*{i};(14,3)*{b+ab-a}
\endxy}=\s_{r=0}^a(-1)^{3a/2-r}x_1^rx_2^{b+ab-r},$$
which implies
$$x_2^{b+ab}=\s_{r=1}^a(-1)^{r+1}x_1^rx_2^{b+ab-r}.$$
We show by downward induction on $1\leq r\leq ab$ that $x_1^rx_2^{b+ab-r}=0$. We can assume $r<b$ so that $ab-r-a\geq 0$. Since
$$0=\ {\fontsize{10}{10}\selectfont\xy
(0,5)*{}; (0,-5)*{} **\crv{(8,0)}?(.1)*{\bu};(5,5)*{}; (5,-5)*{} **\crv{(-3,0)}?(.1)*{\bu};(0,-7)*{i}; (5,-7)*{i};(16.5,3)*{b+ab-r-a};(-2,3)*{r}
\endxy}=\s_{s=0}^a(-1)^{3a/2-s}x_1^{s+r}x_2^{b+ab-r-s},$$
we deduce that
$$x_1^rx_2^{b+ab-r}=\s_{s=1}^a(-1)^{s+1}x_1^{s+r}x_2^{b+ab-r-s}.$$
Note that $s+r<a+b\leq ab+1$, our assertion follows from the induction hypothesis.
\end{proof}
\end{lemma}

\begin{corollary}\label{nil}  Let $\nu\in\N[I]$ with $\tx{ht}(\nu)=n$. Then $x_{k,\ii}$ are nilpotent in $R^\La(\nu)$ for all $\ii\in \Seq(\nu)$ and $1\leq k\leq n$. In particular, $R^\La(\nu)$ is finite dimensional.
\end{corollary}

\begin{remark} 
Note that if $i\in I^0$ and $\La_i=a$, then $R^\La(ni)$ has a basis 
$$\{x_1^{r_1}\cs x_n^{r_n}\tau_\om\mid \om\in S_n,\ 0\leq r_1,\ds,r_n<a\}.$$
\end{remark} 

For each $i\in I$, define the functors
$$\begin{aligned}
&\E_i^\La:R^\La(\nu+i)\Mod\ra R^\La(\nu)\Mod, \ N\ma 1_{\nu,i} N= 1_{\nu,i} R^\La{(\nu+i)}\ot_{R^\La{(\nu+i)}}N, \\
&\F_i^\La:R^\La(\nu)\Mod\ra R^\La(\nu+i)\Mod, \ M\ma R^\La(\nu+i)1_{\nu,i}\ot_{R^\La(\nu)} M.
\end{aligned}$$
Using the nilpotency of $x_k$ (see Corollary \ref{nil}), the following lemma follows by applying the framework of \cite[Section 4]{KK2012} step by step.

\begin{lemma}\label{KK1}
For any $i\in I$, the module $R^\La(\nu+i)1_{\nu,i}$ (resp. $1_{\nu,i}R^\La(\nu+i)$) is projective as a right (resp. left) $R^\La(\nu)$-module. In particular, the functor $\F_i^\La$ is exact, and the functor $\E_i^\La$ sends finitely generated projective modules to finitely generated projective modules.
\end{lemma}

\subsection{Functors $\E^\La_{in},\F^\La_{in}$}\

\vs{1mm}

Let $i\in I^0$ and $n>0$, define the functors
$$\begin{aligned}
& \E_{in}^\La: R^\La(\nu+ni)\Mod\ra R^\La(\nu)\Mod,\  N\ma 1_\nu\ot e_{i,n}\hs{1mm} N,\\
& \F_{in}^\La: R^\La(\nu)\Mod\ra R^\La(\nu+ni)\Mod,\  M\ma (R^\La(\nu+ni)\hs{1mm} 1_\nu\ot e_{i,n})\ot_{R^\La(\nu)}M.\\
\end{aligned}$$
As before, we identify $\E_{in}^\La,\F_{in}^\La$ with their kernels:
$$\E_{in}^\La=1_\nu\ot e_{i,n}\hs{1mm} R^\La(\nu+ni),\quad \F_{in}^\La=R^\La(\nu+ni)\hs{1mm} 1_\nu\ot e_{i,n}.$$

For any $M\in  R^\La(\nu)\Mod$, $N\in R^\La(\nu+ni)\Mod$, we have
$$\HOM_{\nu+ni}((R^\La(\nu+ni)\hs{1mm} 1_\nu\ot e_{i,n})\ot_\nu M,N)\co\HOM_{\nu}(M,\HOM_{\nu+ni}(R^\La(\nu+ni)\hs{1mm} 1_\nu\ot e_{i,n},N)),$$
hence, the functor $\F_{in}^\La$ is left adjoint to $\E_{in}^\La$.

\begin{lemma}\label{TT1}
Let $i\in I^0$. Then $\F_{in}^\La$ (resp. $\E_{in}^\La$) is projective as a right (resp. left) $R^\La(\nu)$-module. In particular, the functor $\F_{in}^\La$ is exact, and $\E_{in}^\La$ takes projective to projective.
\begin{proof}
By Lemma \ref{KK1}, $(\E_i^\La)^n=1_{\nu,ni}R^\La(\nu+ni)$ is a left projective $R^\La(\nu)$-module. Since $\E_{in}^\La$ is a direct summand of $(\E_i^\La)^n$, the assertion follows.
\end{proof}
\end{lemma}

Let $i\in I^0$ and $\nu\in\N[I]$, and assume $\tx{ht}(\nu)=n$. Set $c=(\La-\nu)(h_i)$. A similar argument in \cite[Section 5]{KK2012} gives a natural isomorphism:
$$\E_i^\La\F_i^\La\simeq\F_i^\La\E_i^\La\oplus (\tx{Id}\ot\C[t_i]/(t_i^c)): R^\La(\nu)\Mod\ra R^\La(\nu)\Mod,$$
which is ensured by the following decomposition of $R^\La(\nu)\tx{-}R^\La(\nu)$-bimodules: \begin{equation}\label{tt}
1_{\nu,i}R^\La(\nu+i)1_{\nu,i}=R^\La(\nu)1_{\nu-i,i}\ot_{\nu-i}1_{\nu-i,i}R^\La(\nu)\hs{1mm} \oplus\hs{1mm}  R^\La(\nu)\ot\C[x_{n+1}]/(x_{n+1}^c)
\end{equation}
given by $(x\ot y,z)\ma x\tau_ny+z$.

\begin{proposition}\label{TT2} Let $i\in I^0$ and $\nu\in\N[I]$, and set $c=(\La-\nu)(h_i)$. For any $l,t\geq 1$, we have the following natural isomorphism:
$$\E_{il}^\La\F_{it}^\La\simeq \op_{p=0}^{\min\{l,t\}}\F_{i,t-p}^\La\E_{i,l-p}^\La\ot Z_p^c:\  R(\nu)\Mod\ra R(\nu+(t-l)i)\Mod,$$
where $Z_p^c={\left(\C[x_1,\ds,x_p]/(x_1^c,\ds,x_p^c)\right)}^{S_p}=(P_p^c)^{S_p}$, i.e., the symmetric polynomials in $x_1,\ds,x_p$ in which no $x_k^m$ $(m\geq c)$ appears.
\begin{proof}
Let $\tx{ht}(\nu)=n$. Without loss of generality, assume that $t\geq l$, and set $\nu'=\nu+(t-l)i$.
By iterating the decomposition \eqref{tt}, we obtain
$$\begin{aligned}
& (\E_i^\La)^l(\F_i^\La)^t=1_{\nu',li}R^\La(\nu+ti)1_{\nu,ti}\\
& \phantom{(\E_i^\La)^l(\F_i^\La)^t}=\op_{k=0}^l\op_{\ v\in D_{k,l-k};\ u\in D^\zz_{t-l+k,l-k};\ \om\in S_{l-k}}(R^\La(\nu')\ot v)\x v_k\ot P_{l-k}^c\om \x (R^\La(\nu)\ot u).
\end{aligned}$$

See the diagram below.
$${\fontsize{10}{10}\selectfont\xy
(-6,0)*{}; (-6,12)*{} **\dir{-};(-6,0)*{}; (6,0)*{} **\dir{-};(-6,-9)*{}; (6,-9)*{} **\dir{-};(-6,0)*{}; (-6,-9)*{} **\dir{-};(6,0)*{}; (6,-9)*{} **\dir{-}; (0,-5)*{R^\La(\nu)};
(10,0)*{}; (13,-9)*{} **\dir{-};(13,0)*{}; (19,-9)*{} **\dir{-};(16,0)*{}; (25,-9)*{} **\dir{-};(19,0)*{}; (28,-9)*{} **\dir{-};(22,0)*{}; (10,-9)*{} **\dir{-};(25,0)*{}; (16,-9)*{} **\dir{-};(28,0)*{}; (22,-9)*{} **\dir{-};
(-6,12)*{}; (-6,20)*{} **\dir{-};(9,12)*{}; (9,20)*{} **\dir{-};(-6,12)*{}; (9,12)*{} **\dir{-};(-6,20)*{}; (9,20)*{} **\dir{-};(1.7,15.5)*{R^\La(\nu')};
(13,12)*{}; (16,20)*{} **\dir{-};(16,12)*{}; (22,20)*{} **\dir{-};(19,12)*{}; (28,20)*{} **\dir{-};(22,12)*{}; (13,20)*{} **\dir{-};(25,12)*{}; (19,20)*{} **\dir{-};(28,12)*{}; (25,20)*{} **\dir{-};
(-3,0)*{}; (-3,12)*{} **\dir{-};
(0,0)*{}; (13,12)*{} **\dir{-};
(3,0)*{}; (16,12)*{} **\dir{-};
(6,0)*{}; (19,12)*{} **\dir{-};
(10,0)*{}; (0,12)*{} **\dir{-};
(13,0)*{}; (3,12)*{} **\dir{-};
(16,0)*{}; (6,12)*{} **\dir{-};
(19,0)*{}; (9,12)*{} **\dir{-};
(22,0)*{}; (28,12)*{} **\dir{-}?(0.3)*{\bu};
(25,0)*{}; (25,12)*{} **\crv{(30,6)}?(0.5)*{\bu};
(28,0)*{}; (22,12)*{} **\dir{-}?(0.7)*{\bu};(31,15.5)*{v};(-25,6)*{D_{n,t}\cap D^\zz_{n+t-l,l}\ni v_k};(32,6)*{\om};(31,-4.5)*{u};
(16,10)*{\underbrace{ \quad \  }_{\ \ k}};
\endxy}$$

Let $\wi v=\sum_{v\in D_{k,l-k}}v$ and $\wi u=\sum_{u\in D^\zz_{t-l+k,l-k}}u$. Then we have
$$\begin{aligned}
& \E_{il}^\La\F_{it}^\La=1_{\nu'}\ot e_{i,l}\hs{1mm}R^\La(
\nu+ti)\hs{1mm} 1_{\nu}\ot e_{i,t} \\
& \phantom{\E_{il}^\La\F_{it}^\La}=\op_{k=0}^l\ (R^\La(\nu')\ot e_{i,l})\x v_k\ot P_{l-k}^c \x (R^\La(\nu)\ot e_{i,t})\\
& \phantom{\E_{il}^\La\F_{it}^\La}=\op_{k=0}^l \Big(\big(R^\La(\nu')\hs{1mm} 1_{\nu-ki}\ot e_{i,t-l+k}\big)\ot\wi v\Big) \x v_k\ot e_{i,l-k}P_{l-k}^ce_{i,l-k} \x \Big( \big(1_{\nu-ki}\ot e_{i,k}\hs{1mm} R^\La(\nu)\big)\ot \wi u\Big)\\
& \phantom{\E_{il}^\La\F_{it}^\La}=\op_{k=0}^l \Big(\big(R^\La(\nu')\hs{1mm} 1_{\nu-ki}\ot e_{i,t-l+k}\big)\ot\wi v\Big) \x  v_k\ot Z_{l-k}^c \x \Big( \big(1_{\nu-ki}\ot e_{i,k}\hs{1mm} R^\La(\nu)\big)\ot \wi u\Big).
\end{aligned}$$
For a fixed $k$, the corresponding summand in the above decomposition is the image of
$$\F_{i,t-l+k}^\La\E_{ik}^\La\ot Z_{l-k}^c=\big((R^\La(\nu')\hs{1mm} 1_{\nu-ki}\ot e_{i,t-l+k})\ot_{\nu-ki}(1_{\nu-ki}\ot e_{i,k}\hs{1mm} R^\La(\nu)\big)\ot Z_{l-k}^c$$
under the map
$$(x\ot y,f)\ \ma\ x\ot \wi v\x (v_k\ot f)\x y\ot\wi u,$$
which is injective in the course of the iteration. Consequently, we obtain
$$\E_{il}^\La\F_{it}^\La\simeq\op_{k=0}^l\F_{i,t-l+k}^\La\E_{ik}^\La\ot Z_{l-k}^c\ \ \ov{p=l-k}{\rule{2em}{0.4pt}\kern-2em\raise0.6ex\hbox{\rule{2em}{0.4pt}}}\ \ \op_{k=0}^l\F_{i,t-p}^\La\E_{i,l-p}^\La\ot Z_p^c.$$
\end{proof}
\end{proposition}

\begin{proposition}\label{TT3} We have the following natural isomorphisms:
\begin{itemize}
\item [(1)] For any $(i,n),(j,m)\in\II$, if $a_{ij}=0$, then
$$\E_{in}^\La\E_{jm}^\La\simeq\E_{jm}^\La\E_{in}^\La,\quad \F_{in}^\La\F_{jm}^\La\simeq\F_{jm}^\La\F_{in}^\La.$$
\item [(2)] For $(j,l),(i,t)\in\II$ with $i\neq j$, we have
$$\E_{jl}^\La\F_{it}^\La\simeq q^{-lt a_{ij}}\F_{it}^\La\E_{jl}^\La: R^\La(\nu)\Mod\ra R^\La(\nu+ti-lj)\Mod.$$
\end{itemize}
\begin{proof}
(1) \ If $a_{ij}=0$, then the left multiplication by
$$1_\nu\ot {\fontsize{9}{9}\selectfont\xy
(0,-6)*{}; (12,6)*{} **\dir{-};
(3,-6)*{}; (15,6)*{} **\dir{-};
(6,-6)*{}; (18,6)*{} **\dir{-};
(9,-6)*{}; (0,6)*{} **\dir{-};
(12,-6)*{}; (3,6)*{} **\dir{-};
(15,-6)*{}; (6,6)*{} **\dir{-};
(18,-6)*{}; (9,6)*{} **\dir{-};
(0,-8)*{i}; (3,-8)*{i};(6,-8)*{i};
(9,-8)*{j}; (12,-8)*{j};(15,-8)*{j};(18,-8)*{j};
\endxy}$$
is an isomorphism from
$\E_{in}^\La\E_{jm}^\La=1_{\nu}\ot e_{i,n}\ot e_{j,m}\hs{1mm}R^\La(\nu+ni+mj)$ to
$\E_{jm}^\La\E_{in}^\La=1_{\nu}\ot e_{j,m}\ot e_{i,n}\hs{1mm}R^\La(\nu+ni+mj)$.

(2) Let $\nu'=\nu+ti-lj$. By Proposition \ref{CC1}(3), we have $\E_{jl}\F_{it}\simeq q^{-lt a_{ij}}\F_{it}\E_{jl}$, that is
$$1_{\nu'}\ot e_{j,l}\hs{1mm}R(\nu+ti)\hs{1mm}1_{\nu}\ot e_{i,t}\co q^{-lt a_{ij}}(R(\nu')\hs{1mm}1_{\nu-lj}\ot e_{i,t})\ot_{\nu-lj}(1_{\nu-lj}\ot e_{j,l}\hs{1mm}R(\nu)).$$
Apply the functor $R^\La(\nu')\ot_{R(\nu')}\ \ \bu\ \ {_{R(\nu)}\ot} R^\La(\nu)$ to both sides.
The right-hand side becomes $\F_{it}^\La\E_{jl}^\La$, while the left-hand side is equal to
$$\f{1_{\nu'}\ot e_{j,l}\hs{1mm} R(\nu+ti)\hs{1mm} 1_{\nu}\ot e_{i,t}}{1_{\nu'}\ot e_{j,l}\hs{1mm} R(\nu')a^\La(x_1)R(\nu+ti)\hs{1mm}1_{\nu}\ot e_{i,t}}\ot_{\nu}R^\La(\nu)$$
$$=\f{1_{\nu'}\ot e_{j,l}\hs{1mm}R(\nu+ti)\hs{1mm}1_{\nu}\ot e_{i,t}}{1_{\nu'}\ot e_{j,l}\hs{1mm}R(\nu')a^\La(x_1)R(\nu+ti)\hs{1mm}1_{\nu}\ot e_{i,t}+1_{\nu'}\ot e_{j,l}\hs{1mm}R(\nu+ti)a^\La(x_1)R(\nu)\hs{1mm}1_{\nu}\ot e_{i,t}}$$

Assume $\tx{ht}(\nu)=n$. Note that
$$\begin{aligned}
& \hs{6mm}1_{\nu'}\ot e_{j,l}\hs{1mm}R(\nu+ti)a^\La(x_1)R(\nu+ti)\hs{1mm}1_{\nu}\ot e_{i,t}\\
& =\s_{r=0}^{n+t-1}1_{\nu'}\ot e_{j,l}\hs{1mm} R(\nu+ti)a^\La(x_1)(R(1)\ot R(n+t-1))\x \tau_1\cs\tau_r\hs{1mm}1_{\nu}\ot e_{i,t}\\
& =\s_{r=0}^{n+t-1}1_{\nu'}\ot e_{j,l}\hs{1mm} R(\nu+ti)a^\La(x_1) \tau_1\cs\tau_r\hs{1mm} 1_{\nu}\ot e_{i,t}\\
& =1_{\nu'}\ot e_{j,l}\hs{1mm}R(\nu+ti)a^\La(x_1)R(\nu)\hs{1mm} 1_{\nu}\ot e_{i,t} + 1_{\nu'}\ot e_{j,l}\hs{1mm}R(\nu+ti)a^\La(x_1) \tau_1\cs\tau_n\hs{1mm}1_{\nu}\ot e_{i,t},
\end{aligned}$$
and the second term
$$\begin{aligned}
& \hs{6mm}1_{\nu'}\ot e_{j,l}\hs{1mm}R(\nu+ti)a^\La(x_1) \tau_1\cs\tau_n\hs{1mm}1_{\nu}\ot e_{i,t}\\
& =\s_{r=0}^{n+t-1}1_{\nu'}\ot e_{j,l}\hs{1mm} \tau_r\cs\tau_1\x a^\La(x_1)\ot R(n+t-1) \x\tau_1\cs\tau_n\hs{1mm} 1_{\nu}\ot e_{i,t}\\
& =\s_{r=0}^{n+t-l-1} 1_{\nu'}\ot e_{j,l}\hs{1mm} \tau_r\cs\tau_1\x a^\La(x_1)\ot R(n+t-1) \x\tau_1\cs\tau_n\hs{1mm} 1_{\nu}\ot e_{i,t}\\
& \se\ 1_{\nu'}\ot e_{j,l}\hs{1mm}R(\nu')a^\La(x_1)R(\nu+ti)\hs{1mm}1_{\nu}\ot e_{i,t}.
\end{aligned}$$

$${\fontsize{11}{11}\selectfont
\xy
(0,0)*{}; (15,-11)*{} **\dir{-}; (5,0)*{}; (0,-11)*{} **\dir{-}; (10,0)*{}; (5,-11)*{} **\dir{-}; (15,0)*{}; (10,-11)*{} **\dir{-}; (20,0)*{}; (20,-11)*{} **\dir{-};  (25,0)*{}; (25,-11)*{} **\dir{-};
 (5,-13)*{\underbrace{\qquad \quad }_{\nu}};(20,-14)*{e_{i,t}};
 (0,0)*{}; (0,7)*{} **\dir{-}?(.5)*{\bu}; (5,0)*{}; (5,7)*{} **\dir{-}; (25,0)*{}; (25,7)*{} **\dir{-}; (5,0)*{}; (25,0)*{} **\dir{-}; (5,7)*{}; (25,7)*{} **\dir{-};
 (0,7)*{}; (10,18)*{} **\dir{-};(5,7)*{}; (0,18)*{} **\dir{-};(10,7)*{}; (5,18)*{} **\dir{-};(15,7)*{}; (15,18)*{} **\dir{-};(20,7)*{}; (20,18)*{} **\dir{-};(25,7)*{}; (25,18)*{} **\dir{-};(20,21)*{e_{j,l}};
 (5,20.5)*{\overbrace{\qquad \ \ \ }^{\nu'}}; (15,3)*{\nu +(t-1)i}
\endxy}$$
Hence we conclude that the left-hand side is equal to $\E_{jl}^\La\F_{it}^\La$.
\end{proof}
\end{proposition}

\subsection{Categorification of $V(\La)$}\

\vs{1mm}

Let $G_0(R^\La(\nu))$ (resp. $K_0(R^\La(\nu))$) denote the Grothendieck group of $R^\La(\nu)\fMod$ (resp. $R^\La(\nu)\pMod$). Set
$$R^\La=\op_{\nu\in\N[I]}R^\La(\nu),\quad G_0(R^\La)=\op_{\nu\in\N[I]}G_0(R^\La(\nu)),\quad K_0(R^\La)=\op_{\nu\in\N[I]}K_0(R^\La(\nu)),$$
and $K_0(R^\La)_{\Q(q)}=\Q(q)\ot_{\Z[q,q^{-1}]}K_0(R^\La)$.

For $i\in I^0$, Lemma \ref{TT1} shows that $\E_{in}^\La$ is well defined on $K_0(R^\La)$, and by Proposition \ref{TT2},
$$\E_{il}^\La\F_{it}^\La=\s_{p=0}^{\min\{l,t\}}\Dim Z_p^c\x \F_{i,t-p}^\La\E_{i,l-p}^\La$$
as endomorphisms of $K_0(R^\La(\nu))_{\Q(q)}$, where $c=(\La-\nu)(h_i)$ and $$Z_p^c={\left(\C[x_1,\dots,x_p]/(x_1^c,\ds,x_p^c)\right)}^{S_p}.$$

Note that $Z_p^c$ is determined by all partitions $\lambda$ with $l(\lambda)\leq p$ and $\lambda_1\leq c-1$. We know that the generating function for such partitions is
$$\left\{\begin{matrix} c+p-1 \\ p \end{matrix}\right\}= \f{(1-q^c)(1-q^{c+1})\cs(1-q^{c+p-1})}{(1-q)(1-q^2)\cs(1-q^p)},$$
and therefore
$$\begin{aligned}& \Dim Z_p^c=\f{(1-q^{2c})(1-q^{2(c+1)})\cs(1-q^{2(c+p-1)})}{(1-q^2)(1-q^4)\cs(1-q^{2p})}\\
& \phantom{\Dim Z_p^c}=q^{(c-1)p}\p_{k=1}^p\f{q^{c+k-1}-q^{1-k-c}}{q^k-q^{-k}}:=\bt_p.\end{aligned}$$

\begin{lemma}\label{coeff}
For any $p\geq 1$, $$\bt_p=\ep_p(1-q^{2pc})-\ep_1q^{2c}\bt_{p-1}-\ep_2q^{4c}\bt_{p-2}-\cs-\ep_{p-1}q^{2(p-1)c}\bt_1.$$
\begin{proof}
Using the notations
$$\left\{\begin{matrix} n \\ m \end{matrix}\right\} =\f{(1-q^n)(1-q^{n-1})\cs(1-q^{n-m+1})}{(1-q)(1-q^2)\cs(1-q^m)}$$
for $n\geq m\geq 1$, and $(x;q)_n=(1-x)(1-xq)\cdots (1-xq^{n-1})$ for $n\geq 1$, we have the standard relations
\begin{equation}\label{Gauss}\left\{\begin{matrix} n+1 \\ m \end{matrix}\right\}= q^m\left\{\begin{matrix} n \\ m \end{matrix}\right\}+\left\{\begin{matrix} n \\ m-1 \end{matrix}\right\},\quad (xq^m;q)_{n-m}=\f{(x;q)_n}{(x;q)_m}.\end{equation}

To prove the identity in the lemma, it suffices to show
$$1-q^{pc}=\s_{k=0}^{p-1}q^{kc}\left\{\begin{matrix} p \\ k \end{matrix}\right\}(q^c;q)_{p-k},$$
which can be established easily by induction on $p$, using \eqref{Gauss}.
\end{proof}
\end{lemma}

For $i\in I^-$, a similar argument in \cite[Section 5]{KK2012} leads to
\begin{equation}\label{im}\E_i^\La\F_i^\La=q^{-a_{ii}}\F_i^\La\E_i^\La\oplus \f{1-q^{2c}}{1-q^2}\end{equation}
on $K_0(R^\La(\nu))_{\Q(q)}$, where $c=(\La-\nu)(h_i)$.

For $i\in I^+$, let $c=(\La-\nu)(h_i)$. It was proven in \cite[Section 5]{KK2012} that, on $K_0(R^\La(\nu))_{\Q(q)}$,
\begin{equation}\label{re}\begin{cases}\ \E_i^\La\F_i^\La=q^{-2}\F_i^\La\E_i^\La\oplus \f{1-q^{2c}}{1-q^2}\ & \tx{if}\ c\geq 0, \\ \vs{-4.5mm} \\
\ q^{-2}\F_i^\La\E_i^\La=\E_i^\La\F_i^\La\oplus\f{1-q^{2c}}{q^2-1}
& \tx{if}\ c\leq 0. \end{cases}\end{equation}

Now for any $(i,n)\in\II$ and $h\in P^\vee$, define
$$\FF_{in}^\La=q^{-n(\La-\nu)(h_i)}\F_{in}^\La,\quad q^h=q^{(\La-\nu)(h)}$$
on $K_0(R^\La(\nu))_{\Q(q)}$.

\begin{theorem}\label{TT4} On $K_0(R^\La(\nu))_{\Q(q)}$, we have for any $(j,l),(i,t)\in\II$
$$[\E_{jl}^\La,\FF_{it}^\La]=
\begin{cases}
0  \ \ & \tx{if}\ i\neq j, \\
\sum_{p=1}^{\min\{l,t\}}\ga_p\FF_{i,t-p}^\La\E_{i,l-p}^\La \ \ & \tx{if}\ i=j,
\end{cases}$$
and for $(i,n),(j,m)\in\II$ with $a_{ij}=0$,
$$[\E_{in}^\La,\E_{jm}^\La]=0,\quad [\FF_{in}^\La,\FF_{jm}^\La]=0.$$
In particular, $K_0(R^\La)_{\Q(q)}$ is a weighted $\ovv U$-module ($[P]\in K_0(R^\La(\nu))_{\Q(q)}$ of weight $\La-\nu$).
\begin{proof}
We show the case where $i=j\in I^0$, the other cases follow from Proposition \ref{TT3}, \eqref{im} and \eqref{re}. Recall that
$$\ga_p=\ep_p(K_i^{-p}-K_i^p)-\s_{r=1}^{p-1}\ep_rK_i^r\ga_{p-r}$$
Let $c=(\La-\nu)(h_i)$. Then by Lemma \ref{coeff},
$$\ga_p=\ep_p(q^{-cp}-q^{cp})-\s_{r=1}^{p-1}\ep_rq^{cr}\ga_{p-r}=q^{-cp}\bt_p$$
on $K_0(R^\La(\nu))_{\Q(q)}$. Hence
$$\E_{il}^\La\F_{it}^\La=\s_{p=0}^{\min\{l,t\}}q^{cp}\ga_p\F_{i,t-p}^\La\E_{i,l-p}^\La,$$
and therefore,
$$\E_{il}^\La\FF_{it}^\La=\s_{p=0}^{\min\{l,t\}}q^{c(p-t)}\ga_p \F_{i,t-p}^\La\E_{i,l-p}^\La=\s_{p=0}^{\min\{l,t\}}\ga_p \FF_{i,t-p}^\La\E_{i,l-p}^\La.$$
\end{proof}
\end{theorem}

\begin{proposition}\label{TT5} We have
\begin{itemize}
\item [(1)] If $i\in I^+$, then for any $\nu\in\N[I]$, $R^\La(\nu+ki)=0$ for $k\gg 0$. In particular, $\FF_i^\La$ is locally nilpotent on $K_0(R^\La)_{\Q(q)}$.
\item [(2)] If $i\in I^{\leq 0}$, $\nu\in\N[I]$ with $(\La-\nu)(h_i)=0$, then $R^\La(\nu+ki)=0$ for any $k>0$.
\end{itemize}
In particular, $K_0(R^\La)_{\Q(q)}$ is an integral $\ovv{U}$-module.
\begin{proof}
The proof of (1) is the same as \cite[Lemma 4.3(ii)]{KK2012}. We now prove (2).

Let $i\in I^{\leq 0}$ and suppose $(\La-\nu)(h_i)=0$. Since $a_{ii}\leq 0$ and $a_{ij}\leq 0$ for any $j\neq i$, it follows that $\La(h_i)=0$ and $\nu(h_i)=0$. Hence, by definition, $1_{i,\nu+(k-1)i}R^\La(\nu+ki)=0$. 

For $\ii\in \nu+ki$, let $t$ be the minimal index such that $i_t=i$. An induction on $t$ then shows that  $1_\ii R^\La(\nu+ki)=0$ for any $\ii\in \nu+ki$.
\end{proof}
\end{proposition}

By Lemma \ref{Int}, $K_0(R^\La)_{\Q(q)}$ is an integral $U$-module.
The $\Z[q,q^\zz]$-module $K_0(R^\La)$ and $G_0(R^\La)$ are dual to each other with respect to the Khovanov-Lauda's form
$$K_0(R^\La)\ti G_0(R^\La)\ra \Z[q,q^\zz],\ \ (P,M)\ma\Dim P^\psi\ot_{R^\La}M,$$
and recall that $K_0(R)$ and $G_0(R)$ are similarly dual:
$$K_0(R)\ti G_0(R)\ra \Z[q,q^\zz],\ \ (P,M)\ma\Dim P^\psi\ot_{R}M.$$

The dual map of the canonical embedding
$$G_0(R^\La)\hookrightarrow G_0(R),\ \ [S]\ma [S]$$
is surjective and given by
$$\vp:K_0(R)\twoheadrightarrow K_0(R^\La),\ \ [P]\ma R^\La(\nu)\ot_{R(\nu)}[P].$$

For any $(i,n)\in\II$ and any $[P]\in K_0(R)$, we have
$$\begin{aligned}
& \vp(F_{in}[P])=R^\La(\nu+ni)\ot_{R(\nu+ni)}R(\nu+ni)1_{\nu,ni}\ot_{R(\nu)\ot R(ni)} ([P]\ot R(ni)e_{i,n})\\
& \phantom{\vp(F_{in}[P])}=R^\La(\nu+ni)1_{\nu,ni}\ot_{R(\nu)\ot R(ni)} ([P]\ot R(ni)e_{i,n})\\
& \phantom{\vp(F_{in}[P])}=(R^\La(\nu+ni)\hs{1mm}1_\nu\ot e_{i,n})\ot_{R(\nu)}[P]\\
& \phantom{\vp(F_{in}[P])}=\F_{in}^\La\vp([P])
\end{aligned}$$
It follows that $\vp$ is $U^-$-linear. Therefore, $K_0(R^\La)_{\Q(q)}$ is generated by $[1_\La]$ over $R^\La(0)$ and is an integral highest weight $U$-module of highest weight $\La$. Hence, it is isomorphic to the irreducible highest weight module $V(\La)$.

\begin{theorem}
For each $\La\in P^+$, $K_0(R^\La)_{\Q(q)}\co V(\La)$ as $U$-modules.
\end{theorem}

\vs{7mm}

\bibliographystyle{amsplain}

\end{document}